\documentclass{article}
\usepackage{amssymb,amsmath,bm,geometry,graphics,graphicx,url,color}

\def\no{\noindent}
\def\pmatrix{\left(\begin{array}}
\def\endpmatrix{\end{array}\right)}
\newtheorem{theo}{Theorem}

\newtheorem{defi}{Definition}

\title{Diagonal implicit symplectic
 ERKN methods for solving oscillatory
Hamiltonian systems}%\thanks{This work is supported by Universit\`a di Firenze (project ``Risoluzione numerica di problemi Hamiltoniani ed applicazioni'') and NSFC (Grant No.\,11571128).}}

\author{Mingxue Shi\,
\footnote{School of Mathematical Sciences, Qufu Normal University,
Qufu, Shandong 273165, P.R.China. E-mail:~{\tt 1947785221@qq.com}}
\and Hao Zhang\thanks{School of Mathematical Sciences, Qufu Normal
University, Qufu, Shandong 273165, P.R.China. E-mail:~{\tt
1078836174@qq.com}} \and Bin Wang\thanks{School of Mathematical
Sciences, Qufu Normal University, Qufu, Shandong 273165, P.R.China.
E-mail:~{\tt wangbinmaths@qq.com}} }

\begin{document}
\maketitle

\begin{abstract} This paper studies diagonal implicit symplectic extended
Runge--Kutta--Nystr\"{o}m (ERKN) methods for solving the oscillatory
Hamiltonian system
$H(q,p)=\dfrac{1}{2}p^{T}p+\dfrac{1}{2}q^{T}Mq+U(q)$. Based on
symplectic conditions and order conditions, we construct some
diagonal implicit symplectic ERKN methods. The stability of the
obtained methods is discussed. Three numerical experiments are
carried out and the numerical results
 demonstrate the remarkable numerical behavior of the
 new diagonal implicit symplectic  methods when applied to the oscillatory Hamiltonian
 system.

\medskip
\no{\bf Keywords:} diagonal implicit methods; symplectic methods;
 ERKN methods; oscillatory Hamiltonian systems

\medskip
\no{\bf MSC:} 65L05

\end{abstract}
\section{Introduction}
In this paper, we are concerned with diagonal implicit symplectic
methods for solving the following oscillatory Hamiltonian systems
\begin{equation}
  \dot{q}=\nabla_{p}H(q,p),\qquad
q(t_0)=q_0,\qquad
 \dot{p}=-\nabla_{q}H(q,p),\qquad p(t_0)=p_0\label{a}%
\end{equation}
with the Hamiltonian
\begin{equation}
H(q,p)=\frac{1}{2}p^{T}p+\frac{1}{2}q^{T}Mq+U(q).\label{b}%
\end{equation}
Here $M$ is a $d\times d$ symmetric and positive semi-definite
matrix which implicitly preserves the dominant frequencies of the
system $\eqref{a}$ and function $U(q)$ is a real-valued function
whose second derivatives are continuous. Problem $\eqref{a}$ with
the Hamiltonian of the form $\eqref{b}$ often consists in a wide
variety of applications such as physics, astronomy, molecular
biology,astronomy and classical mechanics (see, e.g.
\cite{cohen2005bit,garcia2002,B.garcia1999,hairer2000,Hairer2002,
Okunbor1994}). Fermi-Pasta-Ulam problem and the spatial
semi-discretization of wave equation with the method of lines are
classcial examples.

 In the past more than twenty
years, some researchers studied the single-frequency problem which
is a kind of special situation of the multi-frequency problem. That
is to say, the multi-frequency problem  $\eqref{a}$  becomes the
single-frequency problem when $M=\omega^{2}$, where $\omega>0$ is
main frequency of the single-frequency problem and may be estimated
in advance. Moreover, the multi-frequency system $\eqref{a}$ is more
complicated than the single-frequency system.  The reasons are as
follows. First, the coefficients of numerical methods for solving
 single-frequency problems depend on $\omega$. However,
$M$ of the multi-frequency system is a $d\times d$ matrix which
implicitly contains many frequencies. Second, symplectic conditions
of methods for solving single-frequency problem can not be extended
to that of methods for solving multi-frequency problem $\eqref{a}$.
So it is important to focus on geometric integration of
multi-frequency problem.

In order to solve the oscillatory Hamiltonian systems $\eqref{a}$,
many novel methods have been developed and studied by many
researchers. For the related work, we can refer the reader to
\cite{wang2016,Wang*2017,Wang*2017*,wu2013-ANM,Wang2017a,wu2015bbb,wu2013-book}.
In \cite{Wu2010}, based on the variation-of-constants formula, Wu et
al. formulate a standard form of  extended Runge--Kutta--Nystr\"{o}m
(ERKN) methods for the oscillatory Hamiltonian system $\eqref{a}$
and present the order conditions via B-series theory. On the other
hand, the importance of geometric numerical integration for the
purpose of preserving some structures of differential equations has
recently become apparent. For the survey of this field, we refer to
the book \cite{Hairer2002}. Ruth \cite{Ruth1983} was the first to
publish the results about canonical numerical methods and proposed
three-stage canonical Runge--Kutta--Nystr\"{o}m (RKN) method of
order three. Some researches following \cite{Ruth1983} about
canonical methods are referred  to
\cite{cohen2005bit,Simos2003,wu-2012-BIT}.  In order to preserve the
symplectic structure of the Hamiltonian system $\eqref{b}$,
symplectic conditions for
 ERKN methods are derived and some novel explicit
 ERKN methods of order up to four are proposed in \cite{wu-2012-BIT}. However, the related work
 about diagonal implicit symplectic  ERKN
methods for solving oscillatory Hamiltonian systems has not been
developed. Therefore, this paper attempts to study  diagonal
implicit symplectic  ERKN methods  which may improve accuracy of the
numerical solutions of the system $\eqref{a}$.

In this paper, we will construct some practical  diagonal implicit
symplectic ERKN methods. It is noted that when matrix
$M\rightarrow\textbf{0}$, these methods will become their
corresponding RKN methods.  Its detailed process is as follows. In
Section \ref{sec:knowledges}, the definition of order conditions and
symplectic conditions of EKRN methods are represented. In Section
\ref{three}, some diagonal implicit symplectic ERKN methods are
constructed. In Section \ref{four}, the stability is analyzed and
 the  obtained methods are compared  with some
 RKN methods by three numerical experiments. In section \ref{five}, we draw the conclusions.

\section{Preliminaries}
\label{sec:knowledges}In order to obtain an effective and practical
numerical scheme for the system $\eqref{a}$, the definition of
 ERKN methods is given in \cite{wu2013-book}.
\begin{defi}
\label{erkn} An s-stage diagonal implicit ERKN method with stepsize
h for solving the
Hamiltonian system \eqref{a} is defined by%
\begin{equation}
\begin{array}
[c]{ll}%
Q_{i} &
=\phi_{0}(c_{i}^{2}V)q_{n}+hc_{i}\phi_{1}(c_{i}^{2}V)p_{n}-h^{2}\textstyle\sum\limits_{j=1}^{i}\bar{a}_{ij}(V)\nabla U(Q_{j}),\ \ \ i=1,\ldots,s,\\
q_{n+1} & =\phi_{0}(V)q_{n}+h\phi_{1}(V)p_{n}-h^{2}\textstyle\sum
\limits_{i=1}^{s}\bar{b}_{i}(V)\nabla U(Q_{i}),\\
p_{n+1} &
=-hM\phi_{1}(V)q_{n}+\phi_{0}(V)p_{n}-h\textstyle\sum\limits_{i=1}^{s}b_{i}(V)\nabla
U(Q_{i}),
\end{array}
  \label{d}%
\end{equation}
where $c_i$ are real constants, $b_{i}(V)$, $\bar{b}_{i}(V),$ and
$\bar{a}_{ij}(V)$ are matrix-valued functions of $V\equiv h^{2}M$,
and
\begin{equation}
\phi_{j}(V):=\sum\limits_{k=0}^{\infty}\dfrac{(-1)^{k}V^{k}}{(2k+j)!},\qquad j=0,1,\ldots i.%
\label{e}%
\end{equation}
The coefficients of $\eqref{d}$ can be displayed in a Butcher
tableau:

\[%
\begin{tabular}
[c]{l}%
\\
\\[2mm]%
\begin{tabular}
[c]{c|c}%
$c$ & $\bar{A}(V)$\\\hline &
$\raisebox{-1.3ex}[0pt]{$\bar{b}^T(V)$}$\\\hline
& $\raisebox{-1.3ex}[0.5pt]{$b^T(V)$}$%
\end{tabular}
$\ \quad=$ $\ $%
\end{tabular}%
\begin{tabular}
[c]{c|ccc}%
$c_{1}$ & $\bar{a}_{11}(V)$ \\
$\vdots$ & $\vdots$ & $\ddots$ \\
$c_{s}$ & $\bar{a}_{s1}(V)$ & $\cdots$ & $\bar{a}_{ss}(V)$\\\hline &
$\raisebox{-1.3ex}[0.5pt]{$\bar{b}_1(V)$}$ &
\raisebox{-1.3ex}{$\cdots$} &
$\raisebox{-1.3ex}[0.5pt]{$\bar{b}_s(V)$}$\\\hline &
$\raisebox{-1.3ex}[1.0pt]{$b_1(V)$}$ &
$\raisebox{-1.3ex}[1.0pt]{$\cdots$}$
& $\raisebox{-1.3ex}[1.0pt]{$b_s(V)$}$%
\end{tabular}
\]
\end{defi}

Based on the new SEN-tree theory, the order conditions of \eqref{d}
can be derived by comparing the series expansions of $q_{n+1}$ and
$p_{n+1}$ in terms of SEN-trees with those of the true solutions
$q(t_{n}+h)$ and $p(t_{n}+h)$, respectively. The following theorem
presents the order conditions of  ERKN methods (see \cite{Wu2010}).
\begin{theo}
\label{abc} The necessary and sufficient conditions for an $s$-stage
 ERKN method $\eqref{d}$ to be
of order $r$ are given by
\begin{equation}\label{f}
\begin{array}
[c]{l}%
\bar{b}^{T}(V)\Phi(\tau)=\dfrac{\rho(\tau)!}{\gamma(\tau)}\phi_{\rho(\tau
)+1}(V)+\mathcal{O}(h^{r-\rho(\tau)}),\ \ \ \rho(\tau)=1,\ldots,r-1,\\
b^{T}(V)\Phi(\tau)=\dfrac{\rho(\tau)!}{\gamma(\tau)}\phi_{\rho(\tau
)}(V)+\mathcal{O}(h^{r+1-\rho(\tau)}),\ \ \ \rho(\tau)=1,\ldots,r,
\end{array}
\end{equation}
where $\tau$ is an extended Nystr\"{o}m tree associated with an
elementary differential $\mathcal{F}(\tau)(q_{n},p_{n})$ of the
function $-\nabla U(q)$ at $q_{n}$.
\end{theo}

The definitions of $\Phi(\tau),\ \rho(\tau),\ \gamma(\tau)$ are
referred to \cite{Wu2010}.

One important property of Hamiltonian systems is that the correspond
flow is symplectic. Therefore, it encourages many researchers to
study symplectic integration which can preserve the symplecticity of
the considered system. The symplectic conditions of ERKN methods are
derived in \cite{wu2013-ANM} and the following theorem states the
symplectic conditions of diagonal implicit ERKN methods for solving
the oscillatory Hamiltonian system \eqref{a}.
\begin{theo}
\label{asdfg} (see \cite{wu-2012-BIT}) An $s$-stage diagonal
implicit ERKN method \eqref{d} is symplectic if its coefficients
satisfy
\begin{equation}
\begin{array}
[c]{ll}
&\phi_{0}(V)b_{i}(V)+V\phi_{1}(V)\bar{b}_{i}(V)=d_{i}\phi_{0}(c_{i}^{2}V),\ \ \ d_{i}\in\mathbb{R},\ \ \ i=1,\ldots,s, \\

&\phi_{0}(V)\bar{b}_{i}(v)+c_{i}d_{i}\phi_{1}(c_{i}^{2}V)=b_{i}(V)\phi_{1}(V),\ \ \ i=1,\ldots,s, \\

&\bar{b}_{j}(V)b_{i}(V)=\bar{b}_{i}(V)b_{j}(V)+d_{i}\bar{a}_{ij}(V),
\ \ \ i=1,\ldots,s,\ \ \ j=1,\ldots,i.
\end{array}
\label{g}%
\end{equation}
\end{theo}
\section{Diagonal implicit symplectic ERKN methods}\label{three}
In this section, we will formulate some diagonal implicit symplectic
ERKN methods with one, two and three stages.

\subsection{One-stage diagonal implicit symplectic methods}\label{ghhgg}
One-stage diagonal implicit ERKN methods can be expressed by a
Butcher tableau:
\[%
\begin{tabular}
[c]{c|ccc}%
$c_{1}$ & $\bar{a}_{11}(V)$ \\\hline &

$\raisebox{-1.3ex}[0.5pt]{$\bar{b}_1(V)$}$\\\hline &

$\raisebox{-1.3ex}[1.0pt]{$b_1(V)$}$%
\end{tabular}
\]
According to the symplectic conditions $\eqref{g}$, it is noted that
this method is symplectic if
\begin{equation}
\begin{array}
[c]{ll}%
\phi_{0}(V)b_{1}(V)+V\phi_{1}(V)\bar{b}_{1}(V)=d_{1}\phi_{0}(c_{1}^{2}V),\\
\phi_{1}(V)b_{1}(V)-\phi_{0}(V)\bar{b}_{1}(V)=c_{1}d_{1}\phi_{1}(c_{1}^{2}V).

\end{array}
  \label{h}%
\end{equation}
 Similarly, according to the order conditions $\eqref{f}$, the
necessary and suffficient conditions for an one-stage
 diagonal implicit ERKN method $\eqref{d}$ to be
of order two are given by
\begin{equation}
\begin{array}
[c]{ll}%
\bar{b}_{1}(V)=\phi_{2}(V)+O(h),\\
b_{1}(V)=\phi_{1}(V)+O(h^{2}),\\
c_{1}b_{1}(V)=\phi_{2}(V)+O(h).
\end{array}
  \label{i}%
\end{equation}
By $\eqref{h}$, we get
\begin{equation*}
\begin{array}
[c]{ll}%
b_{1}(V)=d_{1}\phi_{0}((1-c_{1})^{2}V ),\\
\bar{b}_{1}(V)=d_{1}(1-c_{1})\phi_{1}((1-c_{1})^{2}V).
\end{array}
\label{j}%
\end{equation*}
 Inserting this formula into \eqref{i} yields
\begin{equation*}
\begin{array}
[c]{ll}%
c_{1}=\dfrac{1}{2},\ d_{1}=1.
\end{array}
\label{n}%
\end{equation*}

Due to the second-order order conditions  and  symplectic conditions
do not contain $\bar{a}_{11}(V)$, so the function  $\bar{a}_{11}(V)$
is arbitrary. We consider two kinds of situations.

 \textbf{Case one.}
When $\bar{a}_{11}(V)=\phi_{0}(V)$, the Taylor expansions of other
coefficients   are
\begin{equation}
\begin{aligned}
b_1(V)=&I -\frac{1}{8}V + \frac{1}{384} V^2  - \frac{1}{46080} V^3 +
\cdots,
\\
 \bar{b}_1(V)=&
 \frac{1}{2}I -\frac{1}{48}V + \frac{1}{3840} V^2  -
\frac{1}{645120}
  V^3 + \cdots,
\\
\bar{a}_{11}(V)=& I  -\frac{1}{2}V + \frac{1}{24}V^2
-\frac{1}{720}V^3 +\cdots
.\\
\end{aligned}
\end{equation}
We denote this method by SERKN1s2(1).

 \textbf{Case two.} Choose
$\bar{a}_{11}(V)=\bar{b}_1(V)$ and the Taylor expansions of other
coefficients   are
\begin{equation}
\begin{aligned}
b_1(V)=&I -\frac{1}{8} V + \frac{1}{384} V^2  - \frac{1}{46080} V^3
+ \cdots,
\\
 \bar{b}_1(V)=&
 \frac{1}{2}I -\frac{1}{48}
  V + \frac{1}{3840}V^2 - \frac{1}{645120}
  V^3 + \cdots,\\
\bar{a}_{11}(V)=&\frac{1}{2}I -\frac{1}{48}V + \frac{1}{3840} V^2
-\frac{1}{645120}V^3 +  \cdots.
\\
\end{aligned}
\end{equation}
We denote this method by SERKN1s2(2).

\subsection{Two-stage diagonal implicit  symplectic methods}\label{112233}
We  use a  Butcher  tableau to show two-stage diagonal implicit ERKN
methods:
\[%
\begin{tabular}
[c]{c|ccc}%
$c_{1}$ & $\bar{a}_{11}(V)$ \\
$c_{2}$ & $\bar{a}_{21}(V)$ & $\bar{a}_{22}(V)$\\\hline &
$\raisebox{-1.3ex}[0.5pt]{$\bar{b}_1(V)$}$ &
$\raisebox{-1.3ex}[0.5pt]{$\bar{b}_2(V)$}$\\\hline &
$\raisebox{-1.3ex}[1.0pt]{$b_1(V)$}$ &
$\raisebox{-1.3ex}[1.0pt]{$b_2(V)$}$%
\end{tabular}
\]
By $\eqref{g}$, the symplectic conditions of two-stage diagonal
implicit ERKN methods are given by the following formulas
\begin{equation}
\begin{array}
[c]{ll}
\phi_{0}(V)b_{1}(V)+V\phi_{1}(V)\overline{b}_{1}(V)=d_{1}\phi_{0}(c_{1}^{2}V),\\
\phi_{1}(V)b_{1}(V)-\phi_{0}(V)\overline{b}_{1}(V)=c_{1}d_{1}\phi_{1}(c_{1}^{2}V),\\
\phi_{0}(V)b_{2}(V)+V\phi_{1}(V)\overline{b}_{2}(V)=d_{2}\phi_{0}(c_{2}^{2}V),\\
\phi_{1}(V)b_{2}(V)-\phi_{0}(V)\overline{b}_{2}(V)=c_{2}d_{2}\phi_{1}(c_{2}^{2}V),\\
\overline{b}_{2}(V)b_{1}(V)+d_{2}\overline{a}_{21}(V)=\overline{b}_{1}(V)b_{2}(V).
\end{array}
  \label{q}%
\end{equation}
By $\eqref{f}$, third-order and fourth-order order conditions of
two-stage diagonal implicit ERKN methods respectively are
 \begin{equation}
\begin{array}
[c]{ll} b_{1}(V)+b_{2}(V)=\phi_{1}(V)+O(h^{3}),
\\b_{1}(V)c_{1}+b_{2}(V)c_{2}=\phi_{2}(V)+O(h^{2}),
\\b_{1}(V)c_{1}^{2}+b_{2}(V)c_{2}^{2}=2\phi_{3}(V)+O(h),
\\\overline{b}_{1}(V)+\overline{b}_{2}(V)=\phi_{2}(V)+O(h^{2}),
\\\overline{b}_{1}(V)c_{1}+\overline{b}_{2}(V)c_{2}=\phi_{3}(V)+O(h),
\\b_{1}(V)\overline{a}_{11}(\textbf{0})+b_{2}(V)(\overline{a}_{21}(\textbf{0})+\overline{a}_{22}(\textbf{0}))=\phi_{3}(V)+O(h),
\end{array}
  \label{r}%
\end{equation}
and
\begin{equation}
\begin{array}
[c]{ll} b_{1}(V)+b_{2}(V)=\phi_{1}(V)+O(h^{4}),
\\b_{1}(V)c_{1}+b_{2}(V)c_{2}=\phi_{2}(V)+O(h^{3}),
\\b_{1}(V)c_{1}^{2}+b_{2}(V)c_{2}^{2}=2\phi_{3}(V)+O(h^2),
\\ b_{1}(V)c_{1}^{3}+b_{2}(V)c_{2}^{3}=6\phi_{4}(V)+O(h),
\\\overline{b}_{1}(V)+\overline{b}_{2}(V)=\phi_{2}(V)+O(h^{3}),
\\
\overline{b}_{1}(V)c_{1}+\overline{b}_{2}(V)c_{2}=\phi_{3}(V)+O(h^2),
\\
\overline{b}_{1}(V)c_{1}^2+\overline{b}_{2}(V)c_{2}^2=2\phi_{4}(V)+O(h),\\
b_{1}(V)\overline{a}_{11}(\textbf{0})+b_{2}(V)(\overline{a}_{21}(\textbf{0})+\overline{a}_{22}(\textbf{0}))=\phi_{4}(V)+O(h),
\\
b_{1}(V)\overline{a}_{11}(\textbf{0})+b_{2}(V)(\overline{a}_{21}(\textbf{0})+\overline{a}_{22}(\textbf{0}))=\phi_{3}(V)+O(h^2),
\\
c_{1}b_{1}(V)\overline{a}_{11}(\textbf{0})+c_{2}b_{2}(V)(\overline{a}_{21}(\textbf{0})+\overline{a}_{22}(\textbf{0}))=3\phi_{4}(V)+O(h),
\\
\end{array}
\label{s}%
\end{equation}
\begin{equation*}
\begin{array}
[c]{ll}
c_{1}b_{1}(V)\overline{a}_{11}(\textbf{0})+b_{2}(V)(c_{1}\overline{a}_{21}(\textbf{0})+c_{2}\overline{a}_{22}(\textbf{0}))=\phi_{4}(V)+O(h).
\end{array}
\end{equation*}
By the first four formulas of $\eqref{q}$, we obtain
\begin{equation}
\begin{array}
[c]{ll}
\bar{b}_{1}(V)=b_{1}(1-c_{1})\dfrac{\phi_{1}((1-c_{1})^{2}V)}{\phi_{0}(1-c_{1})^{2}V},\
\
&\bar{b}_{2}(V)=b_{2}(1-c_{2})\dfrac{\phi_{1}((1-c_{2})^{2}V)}{\phi_{0}(1-c_{2})^{2}V},\\
b_{1}(V)=d_{1}\phi_{0}((-1+c_{1})^{2}V),\ \
&b_{2}(V)=d_{2}\phi_{0}((-1+c_{2})^{2}V).
\end{array}
\label{t}%
\end{equation}
By the last formula of  $\eqref{q}$, we have
\begin{equation*}
\begin{array}
[c]{ll}
\bar{a}_{21}(V)=\dfrac{b_{2}(V)\bar{b}_{1}(V)-b_{1}(V)\bar{b}_{2}(V)}{d_{2}}.
\end{array}
\label{u}%
\end{equation*}

\textbf{Case one.} Inserting these formulas of $\eqref{t}$ into the
first five formulas of $\eqref{r}$ yields
\begin{equation*}
\begin{array}
[c]{ll} d_{1}=\dfrac{1-2c_{2}}{2(c_{1}-c_{2})},\ \
 d_{2}=\dfrac{-1+2c_{1}}{2(c_{1}-c_{2})},\ \
 c_{2}=\dfrac{2-3c_{1}}{3-6c_{1}},
\end{array}
\label{x}%
\end{equation*}
where  $c_{1}$ is a parameter.

 Considering
$\bar{a}_{22}(V)=\bar{a}_{11}(V)$ and the following formula (which
is obtained from the last formula of
 $\eqref{r}$)
$$b_{1}(V)\overline{a}_{11}(V)+b_{2}(\overline{a}_{21}(V)+\overline{a}_{22}(V))=\phi_{3}(V),$$
we obtain
$$\bar{a}_{22}(V)=\bar{a}_{11}(V)=\dfrac{-\bar{a}_{21}(V)b_{2}(V)+\phi_{3}(V)}{b_{1}(V)+b_{2}(V)}.$$
We choose $c_{1}=\dfrac{1}{5}$ and then get $c_2=\dfrac{7}{9},$\ \
 $d_1=\dfrac{25}{52},$\ \ $d_2=\dfrac{27}{52}$. The Taylor expansions
of $b_{i}(V),$\ \ $\bar{b}_{i}(V),$\ \ $\bar{a}_{ij}(V)$ are
\begin{equation}
\begin{aligned}
b_1(V)=&\frac{25}{52}I -\frac{2}{13} V + \frac{8}{975} V^2 -
\frac{64}{365625} V^3 +  \cdots,
\\
b_2(V)=&\frac{27}{52}I -\frac{1}{78} V + \frac{1}{18954} V^2
-\frac{1}{11514555}V^{3} +\cdots,
\\
 \bar{b}_1(V)=&
 \frac{5}{13}I  -\frac{8}{195}
  V + \frac{32}{24375} V^2 - \frac{256}{12796875} V^3 +\cdots
,\\
\bar{b}_2(V)=&
 \frac{3}{26}I -\frac{1}{1053}
  V +  \frac{1}{426465} V^2  - \frac{2}{725416965} V^3 +\cdots
,\\
\bar{a}_{11}(V)=&
 \frac{7}{312}I +\frac{589}{84240}V +
 \frac{9931609}{11941020000}V^2+
 \frac{25321869691}{290166786000000}V^3 +\cdots
,\\
\bar{a}_{21}(V)=& \frac{5}{18}I -\frac{169}{10935}V +
\frac{28561}{110716875}V^2 -\frac{9653618}{4708235109375}V^3  +
\cdots
,\\
\bar{a}_{22}(V)=&
  \frac{7}{312}I +\frac{589}{84240}V
+ \frac{9931609}{11941020000}V^2 +
\frac{25321869691}{290166786000000}V^3 + \cdots
.\\
\end{aligned}
\end{equation}
We denote this method by SERKN2s3.

\textbf{Case two.} Inserting  $\eqref{t}$ into the first seven
formulas of $\eqref{s}$ yields
\begin{equation*}
\begin{array}
[c]{ll} d_{1}=\dfrac{1-2c_{2}}{2(c_{1}-c_{2})},\ \
 d_{2}=\dfrac{-1+2c_{1}}{2(c_{1}-c_{2})}.
\end{array}
\label{ppmn}%
\end{equation*}
Considering $c_{1}=\dfrac{3-\sqrt{3}}{6}$,
$c_{2}=\dfrac{3+\sqrt{3}}{6}$, we get $d_{1}=\dfrac{1}{2}$,
$d_{2}=\dfrac{1}{2}$. According to the eighth and ninth formula of
\eqref{s}, the follow results are obtained:
\begin{equation*}
\begin{array}
[c]{ll}
\bar{a}_{11}(V)=\dfrac{\bar{b}_{2}(V)\phi_{3}(V)-b_{2}(V)\phi_{4}(V)}{-b_{2}(V)\bar{b}_{1}(V)+b_{1}(V)\bar{b}_{2}(V)},\\
\bar{a}_{22}(V)=\dfrac{\bar{a}_{21}(V)b_{2}(V)\bar{b}_{1}(V)-\bar{a}_{21}(V)b_{1}(V)\bar{b}_{2}(V)-
\bar{b}_{1}(V)\phi_{3}(V)+b_{1}(V)\phi_{4}(V)}{-b_{2}(V)\bar{b}_{1}(V)+b_{1}(V)\bar{b}_{2}(V)}.
\end{array}
\label{ppmn}%
\end{equation*}
The Taylor expansions of other coefficients are
\begin{equation}
\begin{aligned}
b_1(V)=&\frac{1}{2}I -\frac{2+\sqrt{3}}{24} V +
\frac{7+4\sqrt{3}}{1728} V^2 -\frac{26+15\sqrt{3}}{311040} V^3 +
\cdots
,\\
b_2(V)=&\frac{1}{2}I + \frac{-2+\sqrt{3}}{24}V +
\frac{7-4\sqrt{3}}{1728} V^2 + \frac{-26+15\sqrt{3}}{311040} V^3
+\cdots ,
\\
 \bar{b}_1(V)=&
 \frac{3+\sqrt{3}}{12}I -\frac{9+5\sqrt{3}}{432}
  V + \frac{33+19\sqrt{3}}{51840}V^2 -\frac{123+71\sqrt{3}}{13063680}
  V^3  +\cdots
,\\
 \bar{b}_2(V)=&
 \frac{3-\sqrt{3}}{12}I  +\frac{-9+5\sqrt{3}}{432}
  V + \frac{33-19\sqrt{3}}{51840} V^2 + \frac{-123+71\sqrt{3}}{13063680}
  V^3 +\cdots
,\\
\bar{a}_{11}(V)=& \frac{2-\sqrt{3}}{12}I +\frac{-3+2\sqrt{3}}{2160}V
+ \frac{6-\sqrt{3}}{1088640}V^2  +\frac{15+2\sqrt{3}}{195955200}V^3
+\cdots
,\\
\bar{a}_{21}(V)=& \frac{1}{2\sqrt{3}}I -\frac{1}{36\sqrt{3}}V +
\frac{1}{2160\sqrt{3}}V^2 -\frac{1}{272160\sqrt{3}}V^3  +\cdots
,\\
\bar{a}_{22}(V)=& \frac{2-\sqrt{3}}{12}I +\frac{-1+6\sqrt{3}}{720}V
+ \frac{6-167\sqrt{3}}{1088640} V^2  +
\frac{15+238\sqrt{3}}{195955200}V^3 +\cdots
.\\
\end{aligned}
\end{equation}
We denote this method by SERKN2s4 and this method is proved to
satisfy all the order conditions and symplectic conditions.

\subsection{Three-stage diagonal implicit  symplectic methods}\label{334455}
The following Butcher tableau is given to describe three-stage
diagonal implicit methods:
\[%
\begin{tabular}
[c]{c|ccc}%
$c_{1}$ & $\bar{a}_{11}(V)$ \\
$c_{2}$ & $\bar{a}_{21}(V)$ & $\bar{a}_{22}(V)$ \\
$c_{3}$ & $\bar{a}_{31}(V)$ & $\bar{a}_{32}(V)$ &
$\bar{a}_{33}(V)$\\\hline &
$\raisebox{-1.3ex}[0.5pt]{$\bar{b}_1(V)$}$ &
\raisebox{-1.3ex}{$\bar{b}_2(V)$} &
$\raisebox{-1.3ex}[0.5pt]{$\bar{b}_3(V)$}$\\\hline &
$\raisebox{-1.3ex}[1.0pt]{$b_1(V)$}$ &
$\raisebox{-1.3ex}[1.0pt]{$b_2(V)$}$
& $\raisebox{-1.3ex}[1.0pt]{$b_3(V)$}$%
\end{tabular}
\]
By $\eqref{g}$ and $\eqref{f}$, the symplectic condition and
fourth-order order conditions of three-stage diagonal implicit ERKN
methods respectively are given by the following formulas
\begin{equation}
\begin{array}
[c]{ll}
\phi_{0}(V)b_{1}(V)+V\phi_{1}(V)\overline{b}_{1}(V)=d_{1}\phi_{0}(c_{1}^{2}V),
 \\\phi_{1}(V)b_{1}(V)-\phi_{0}(V)\overline{b}_{1}(V)=c_{1}d_{1}\phi_{1}(c_{1}^{2}V),
 \\\phi_{0}(V)b_{2}(V)+V\phi_{1}(V)\overline{b}_{2}(V)=d_{2}\phi_{0}(c_{2}^{2}V),
 \\\phi_{1}(V)b_{2}(V)-\phi_{0}\overline{b}_{2}(V)=c_{2}d_{2}\phi_{1}(c_{2}^{2}V),
 \\\phi_{0}(V)b_{3}(V)+V\phi_{1}\overline{b}_{3}(V)=d_{3}\phi_{0}(c_{3}^{2}V),
 \\\phi_{1}(V)b_{3}(V)-\phi_{0}\overline{b}_{3}(V)=c_{3}d_{3}\phi_{1}(c_{3}^{2}V),
 \\\overline{b}_{1}(V)b_{2}(V)=\overline{b}_{2}(V)b_{1}(V)+d_{2}\overline{a}_{21}(V),
 \\\overline{b}_{1}(V)b_{3}(V)=\overline{b}_{3}(V)b_{1}(V)+d_{3}\overline{a}_{31}(V),
 \\\overline{b}_{2}(V)b_{3}(V)=\overline{b}_{3}(V)b_{2}(V)+d_{3}\overline{a}_{32}(V),
\end{array}
 \label{erkn123}%
\end{equation}
and
\begin{equation}
\begin{array}
[c]{ll}b_{1}(V)+b_{2}(V)+b_{3}(V)=\phi_{1}(V)+O(h^{4}),
\\b_{1}(V)c_{1}+b_{2}(V)c_{2}+b_{3}(V)c_{3}=\phi_{2}(V)+O(h^{3}),
\\b_{1}(V)c_{1}^{2}+b_{2}(V)c_{2}^{2}+b_{3}(V)c_{3}^{2}=2\phi_{3}(V)+O(h^{2}),
\\b_{1}(V)c_{1}^{3}+b_{2}(V)c_{2}^{3}+b_{3}(V)c_{3}^{3}=6\phi_{4}(V)+O(h),
\\\overline{b}_{1}(V)+\overline{b}_{2}(V)+\overline{b}_{3}(V)=\phi_{2}(V)+O(h^{3}),
\\\overline{b}_{1}(V)c_{1}+\overline{b}_{2}(V)c_{2}+\overline{b}_{3}(V)c_{3}=\phi_{3}(V)+O(h^{2}),
\\\overline{b}_{1}(V)c_{1}^{2}+\overline{b}_{2}(V)c_{2}^{2}+\overline{b}_{3}(V)c_{3}^{2}=2\phi_{4}(V)+O(h),\\
\overline{b}_{1}(V)\overline{a}_{11}(\textbf{0})+\overline{b}_{2}(V)(\overline{a}_{21}(\textbf{0})+\overline{a}_{22}(\textbf{0}))+\overline{b}_{3}(V)(\overline{a}_{31}(\textbf{0})+\overline{a}_{32}(\textbf{0})+\overline{a}_{33}(\textbf{0}))=\phi_{4}(V)+O(h),
\\b_{1}(V)\overline{a}_{11}(\textbf{0})+b_{2}(V)(\overline{a}_{21}(\textbf{0})+\overline{a}_{22}(\textbf{0}))+b_{3}(V)(\overline{a}_{31}(\textbf{0})+\overline{a}_{32}(\textbf{0})+\overline{a}_{33}(\textbf{0}))=\phi_{3}(V)+O(h^{2}),\\
b_{1}(V)c_{1}\overline{a}_{11}(\textbf{0})+b_{2}(V)c_{2}(\overline{a}_{21}(\textbf{0})+\overline{a}_{22}(\textbf{0}))+b_{3}(V)c_{3}(\overline{a}_{31}(\textbf{0})+\overline{a}_{32}(\textbf{0})+\overline{a}_{33}(\textbf{0}))=3\phi_{4}(V)+O(h),
\\
\end{array}
 \label{Schems100}%
\end{equation}
\begin{equation*}
\begin{array}
[c]{ll}
b_{1}(V)c_{1}\overline{a}_{11}(\textbf{0})+b_{2}(V)(c_{1}\overline{a}_{21}(\textbf{0})+c_{2}\overline{a}_{22}(\textbf{0}))+b_{3}(V)(c_{1}\overline{a}_{31}(\textbf{0})+c_{2}\overline{a}_{32}(\textbf{0})+c_{3}\overline{a}_{33}(\textbf{0}))=\phi_{4}(V)
+O(h).
\end{array}
\end{equation*}
On the one hand, $b_{i}$ and $\bar{b}_{i},$\ \ $i=1,2,3$ can be
obtained  by solving the first six formulas of $\eqref{erkn123}$ as
follows
\begin{equation}
\begin{array}
[c]{ll}
\bar{b}_{1}(V)=\dfrac{b_{1}(\phi_{0}(c_{1}^2V)\phi_{1}(V)-c_{1}\phi_{0}(V)\phi_{1}(c_{1}^2V))}{\phi_{0}(V)\phi_{0}(c_{1}^2V)+c_{1}V\phi_{1}(V)\phi_{1}(c_{1}^2V)},\\
\bar{b}_{2}(V)=\dfrac{b_{2}(\phi_{0}(c_{2}^2V)\phi_{1}(V)-c_{2}\phi_{0}(V)\phi_{1}(c_{2}^2V))}{\phi_{0}(V)\phi_{0}(c_{2}^2V)+c_{2}V\phi_{1}(V)\phi_{1}(c_{2}^2V)},\\
\bar{b}_{3}(V)=\dfrac{b_{3}(\phi_{0}(c_{3}^2V)\phi_{1}(V)-c_{3}\phi_{0}(V)\phi_{1}(c_{3}^2V))}{\phi_{0}(V)\phi_{0}(c_{3}^2V)+c_{3}V\phi_{1}(V)\phi_{1}(c_{3}^2V)},\\
b_{1}(V)=\dfrac{d_{1}(\phi_{0}(V)\phi_{0}(c_{1}^2V)+c_{1}V\phi_{1}(V)\phi_{1}(c_{1}^2V))}{\phi_{0}(V)^2+V\phi_{1}(V)^2},\\
b_{2}(V)=\dfrac{d_{2}(\phi_{0}(V)\phi_{0}(c_{2}^2V)+c_{2}V\phi_{1}(V)\phi_{1}(c_{2}^2V))}{\phi_{0}(V)^2+V\phi_{1}(V)^2},\\
b_{3}(V)=\dfrac{d_{3}(\phi_{0}(V)\phi_{0}(c_{3}^2V)+c_{3}V\phi_{1}(V)\phi_{1}(c_{3}^2V))}{\phi_{0}(V)^2+V\phi_{1}(V)^2}.
\end{array}
  \label{scvg}%
\end{equation}
 On the other hand, by the last three formulas of
$\eqref{erkn123}$, we   obtain
\begin{equation}
\begin{array}
[c]{ll}
\bar{a}_{21}(V)=\dfrac{b_{2}(V)\bar{b}_{1}(V)-b_{1}(V)\bar{b}_{2}(V)}{d_{2}},\\
\bar{a}_{31}(V)=\dfrac{b_{3}(V)\bar{b}_{1}(V)-b_{1}(V)\bar{b}_{3}(V)}{d_{3}},\\
\bar{a}_{32}(V)=\dfrac{b_{3}(V)\bar{b}_{2}(V)-b_{2}(V)\bar{b}_{3}(V)}{d_{3}}.
\end{array}
  \label{ekrnA}%
\end{equation}
In the above formulas, $c_{i}$ and $d_{i},\ \ i=1,2,3 $ are
parameters.

 It is noted that $\bar{a}_{ii}$\ \ $(i=1,2,3)$ can be obtained by solving the eighth,
 ninth, tenth formulas of $\eqref{Schems100}$ as follows
\begin{equation}
\begin{array}
[c]{ll}
\bar{a}_{11}(V)=\dfrac{((-b_{2}(V)\bar{b}_{3}(V)c_{2}+b_{3}(V)\bar{b}_{2}(V)c_{3})\phi_{3}(V)}{b_{2}(V)b_{3}(V)\bar{b}_{1}(V)(c_{2}-c_{3})+b_{1}(V)(b_{2}(V)\bar{b}_{3}(V)(c_{1}-c_{2})+b_{3}(V)\bar{b}_{2}(V)(-c_{1}+c_{3}))}\\
\ \ \ \ \ \ \ \ \ \ \ \ \ \ \
+\dfrac{(3b_{2}(V)\bar{b}_{3}(V)+b_{3}(V)(-3\bar{b}_{2}(V)+b_{2}(V)c_{2}(V)
-b_{2}(V)c_{3}(V)))\phi_{4}(V))}{b_{2}(V)b_{3}(V)\bar{b}_{1}(V)(c_{2}-c_{3})+b_{1}(V)(b_{2}(V)\bar{b}_{3}(V)(c_{1}-c_{2})+b_{3}(V)\bar{b}_{2}(V)(-c_{1}+c_{3}))},\\
\bar{a}_{22}(V)=\dfrac{\bar{a}_{21}(V)(b_{1}(V)(b_{2}(V)\bar{b}_{3}(V)(-c_{1}+c_{2})+b_{3}(V)\bar{b}_{2}(V)(c_{1}-c_{3})}{b_{2}(V)b_{3}(V)\bar{b}_{1}(c_{2}-c_{3})
+b_{1}(V)(b_{2}(V)\bar{b}_{3}(V)(c_{1}-c_{2})+b_{3}(V)\bar{b}_{2}(V)(-c_{1}+c_{3}))}\\
\ \ \ \ \ \ \ \ \ \ \ \ \ \ \
+\dfrac{b_{2}(V)b_{3}(V)\bar{b}_{1}(V)(-c_{2}+c_{3}))
+(b_{1}(V)\bar{b}_{3}(V)c_{1}-b_{3}(V)\bar{b}_{1}(V)c_{3})\phi_{3}(V)}{b_{2}(V)b_{3}(V)\bar{b}_{1}(V)(c_{2}-c_{3})
+b_{1}(V)(b_{2}(V)\bar{b}_{3}(V)(c_{1}-c_{2})+b_{3}(V)\bar{b}_{2}(V)(-c_{1}+c_{3}))}\\
\ \ \ \ \ \ \ \ \ \ \ \ \ \ \
+\dfrac{(-3b_{1}(V)\bar{b}_{3}(V)+b_{3}(V)(3\bar{b}_{1}(V)-b_{1}(V)c_{1}+b_{1}(V)c_{3}))\phi_{4}}{b_{2}(V)b_{3}(V)\bar{b}_{1}(V)(c_{2}-c_{3})
+b_{1}(V)(b_{2}(V)\bar{b}_{3}(V)(c_{1}-c_{2})+b_{3}(V)\bar{b}_{2}(V)(-c_{1}+c_{3}))}
 ,\\

\bar{a}_{33}(V)=\dfrac{(\bar{a}_{31}(V)+a_{32}(V))(b_{1}(V)(b_{2}(V)\bar{b}_{3}(V)(-c_{1}+c_{2})+b_{3}(V)\bar{b}_{2}(V)(c_{1}-c_{3})}{b_{2}(V)b_{3}(V)\bar{b}_{1}(V)(c_{2}-c_{3})+b_{1}(V)(b_{2}(V)\bar{b}_{3}(V)(c_{1}-c_{2})+b_{3}(V)\bar{b}_{2}(V)(-c_{1}+c_{3}))}\\
\ \ \ \ \ \ \ \ \ \ \ \ \ \ \ +\dfrac{b_{2}(V)b_{3}(V)\bar{b}_{1}(V)
\times(-c_{2}+c_{3}))+(-b_{1}(V)\bar{b}_{2}(V)c_{1}+b_{2}(V)\bar{b}_{1}(V)c_{2})\phi_{3}(V)}{b_{2}(V)b_{3}(V)\bar{b}_{1}(V)(c_{2}-c_{3})+b_{1}(V)(b_{2}(V)\bar{b}_{3}(V)(c_{1}-c_{2})+b_{3}(V)\bar{b}_{2}(V)(-c_{1}+c_{3}))}\\
\ \ \ \ \ \ \ \ \ \ \ \ \ \ \
+\dfrac{3b_{1}(V)\bar{b}_{2}(V)+b_{2}(V)(-3\bar{b}_{1}(V)+b_{1}(V)c_{1}-b_{1}(V)c_{2}))\phi_{4}(V)}{b_{2}(V)b_{3}(V)\bar{b}_{1}(V)(c_{2}-c_{3})+b_{1}(V)(b_{2}(V)\bar{b}_{3}(V)(c_{1}-c_{2})+b_{3}(V)\bar{b}_{2}(V)(-c_{1}+c_{3}))}.
\end{array}
\label{ekrnB}%
\end{equation}
In order to fulfill  the first seven formulas of $\eqref{Schems100}$
, the following conditions should be satisfied:
\begin{equation}
\begin{array}
[c]{ll}
d_{1}=\dfrac{2-3c_{3}+c_{2}(-3+6c_{3})}{6(c_{1}-c_{2})(c_{1}-c_{3})},\\
d_{2}=\dfrac{-2+c_{1}(3-6c_{3})+3c_{3}}{6(c_{1}-c_{2})(c_{2}-c_{3})},\\
d_{3}=\dfrac{-2+c_{1}(3-6c_{2})+3c_{2}}{6(c_{1}-c_{3})(-c_{2}+c_{3})},\\
c_{3}=\dfrac{3-4c_{1}-4c_{2}+6c_{1}c_{2}}{4-6c_{1}-6c_{2}+12c_{1}c_{2}},\\
\end{array}
  \label{ekrnC}%
\end{equation}
where $c_{1} $ and $c_{2}$ are parameters.

\textbf{Case one.}  We choose
\begin{equation}
\begin{array}
[c]{ll}
c_{1}=\dfrac{5-\sqrt{15}}{10},\ \ c_{2}=\dfrac{1}{2},\\
 \end{array}
  \label{ekrnD}%
\end{equation}
and then get $c_{3}=\dfrac{5+\sqrt{15}}{10}$. Under the conditions
of $\eqref{ekrnD}$, the Taylor series expansions of $b_{i}(V)$,\ \
$\bar{b}_{i}(V)$,\ \ $\bar{a}_{ij}(V)$\ \ $(i=1,2,3,$\ \ $j=1\ldots
i)$ respectively are
\begin{equation}
\begin{aligned}
 b_1(V)=&
\frac{5}{18}I - \frac{4+\sqrt{15}}{72} V +
\frac{31+8\sqrt{15}}{8640} V^2 -\frac{244+63\sqrt{15}}{2592000} V^3
+ \cdots,
\\
b_2(V)=&\frac{4}{9}I -\frac{1}{18}V  + \frac{1}{864} V^2
-\frac{1}{103680} V^3   +\cdots,
\\
b_3(V)=& \frac{5}{18}I  + \frac{-4+\sqrt{15}}{72} V +
\frac{31-8\sqrt{15}}{8640}V^2 +\frac{-244+63\sqrt{15}}{2592000} V^3  +\cdots,\\
 \bar{b}_1(V)=&
 \frac{5+\sqrt{15}}{36}I - \frac{35+9\sqrt{15}}{2160}
  V + \frac{275+71\sqrt{15}}{432000} V^2 - \frac{2165+559\sqrt{15}}{181440000} V^3  +\cdots
,\\
 \bar{b}_2(V)=&
 \frac{2}{9}I -\frac{1}{108}
  V + \frac{1}{8640} V^2 - \frac{1}{1451520} V^3  +\cdots
,\\
 \bar{b}_3(V)=&
 \frac{5-\sqrt{15}}{36}I + \frac{-35+9\sqrt{15}}{2160}
  V + \frac{275-71\sqrt{15}}{432000} V^2 + \frac{-2165+559\sqrt{15}}{181440000} V^3  +\cdots
,\\
\bar{a}_{11}(V)=& \frac{4-\sqrt{15}}{20}I -\frac{3}{14000}V  -
\frac{28+9\sqrt{15}}{5040000} V^2
-\frac{5871+440\sqrt{15}}{16632000000}V^3 +\cdots
,\\
\bar{a}_{21}(V)=& \frac{\sqrt{15}}{36}I -\frac{1}{96\sqrt{15}}V
+\frac{1}{12800\sqrt{15}}V^2 -\frac{1}{3584000\sqrt{15}}V^3   +
\cdots
,\\
\bar{a}_{22}(V)=& \frac{9-2\sqrt{15}}{72}I
+\frac{3}{11200}V+\frac{1}{96\sqrt{15}}V +
\frac{29-42\sqrt{15}}{8064000}V^2 +
\frac{-3667+330\sqrt{15}}{17740800000}V^3  +\cdots
,\\
\bar{a}_{31}(V)=& \frac{\sqrt{15}}{18}I -\frac{1}{12\sqrt{15}}V
+\frac{1}{400\sqrt{15}}V^2 -\frac{1}{28000\sqrt{15}}V^3
 +\cdots
,\\
\bar{a}_{32}(V)=& \frac{2}{3\sqrt{15}}I -\frac{1}{60\sqrt{15}}V
+\frac{1}{8000\sqrt{15}}V^2 -\frac{1}{2240000\sqrt{15}}V^3 +\cdots
,\\
\bar{a}_{33}(V)=& \frac{4-\sqrt{15}}{20}I
+\frac{-9+280\sqrt{15}}{42000}V -\frac{28+873\sqrt{15}}{5040000}V^2
+\frac{-5871+40535\sqrt{15}}{16632000000}V^3  +\cdots
.\\
\end{aligned}
\end{equation}
We denote this method by SERKN3s4(1), and this method is proved to
satisfy all the order conditions and symplectic conditions.

 \textbf{Case two.} Choose
\begin{equation}
\begin{array}
[c]{ll}c_{1}=\dfrac{5+\sqrt{15}}{10},\ \
c_{2}=\dfrac{5-\sqrt{15}}{10},
 \end{array}
  \label{ekrnE}%
\end{equation}
and then we get $c_{3}=\dfrac{1}{2}$. Under the conditions of
$\eqref{ekrnE}$,
 the Taylor series expansions of
other coefficients respectively are
\begin{equation}
\begin{aligned}
b_1(V)=&   \frac{5}{18}I  + \frac{-4+\sqrt{15}}{72} V +
\frac{31-8\sqrt{15}}{8640} V^2 + \frac{-244+63\sqrt{15}}{2592000}
V^3 +\cdots  ,
\\
b_2(V)=& \frac{5}{18}I - \frac{4+\sqrt{15}}{72} V  +
\frac{31+8\sqrt{15}}{8640}V^2 -\frac{244+63\sqrt{15}}{2592000} V^3
+\cdots,
\\
b_3(V)=&\frac{4}{9}I-\frac{1}{18} V +\frac{1}{864}
V^2-\frac{1}{103680} V^3 +\cdots
,\\
\bar{b}_1(V)=& \frac{5-\sqrt{15}}{36}I + \frac{-35+9\sqrt{15}}{2160}
  V +
\frac{275-71\sqrt{15}}{432000} V^2 +
\frac{-2165+559\sqrt{15}}{181440000}V^3  +\cdots
,\\
\bar{b}_2(V)=& \frac{5+\sqrt{15}}{36}I - \frac{35+9\sqrt{15}}{2160}
  V +
\frac{275+71\sqrt{15}}{432000} V^2
-\frac{2165+559\sqrt{15}}{181440000} V^3  +\cdots
,\\
\end{aligned}
\end{equation}
\begin{equation}
\begin{aligned}
\bar{b}_3(V)=&
 \frac{2}{9}I -\frac{1}{108}
  V +\frac{1}{8640} V^2 -\frac{1}{1451520} V^3  +\cdots
,\\
\bar{a}_{11}(V)=& \frac{4+\sqrt{15}}{20}I -\frac{3}{14000}V +
\frac{-28+9\sqrt{15}}{5040000} V^2
+\frac{-5871+440\sqrt{15}}{16632000000}V^3  +\cdots
,\\
\bar{a}_{21}(V)=&
 -\frac{\sqrt{15}}{18}I
+\frac{1}{12\sqrt{15}}V -\frac{1}{400\sqrt{15}}V^2 +
\frac{1}{28000\sqrt{15}}V^3 +\cdots
,\\
\bar{a}_{22}(V)=& \frac{36+\sqrt{15}}{180}I -\frac{1}{12\sqrt{15}}V
-\frac{3}{14000}V + \frac{-28+831\sqrt{15}}{5040000} V^2
-\frac{5871+40040\sqrt{15}}{16632000000}V^3 +\cdots
,\\
\bar{a}_{31}(V)=&
 -\frac{\sqrt{15}}{36}I
+\frac{1}{96\sqrt{15}}V -\frac{1}{12800\sqrt{15}}V^2 +
\frac{1}{3584000\sqrt{15}}V^3 +\cdots
,\\
\bar{a}_{32}(V)=&\frac{\sqrt{15}}{36}I -\frac{1}{96\sqrt{15}}V +
\frac{1}{12800\sqrt{15}}V^2 -\frac{1}{3584000\sqrt{15}}V^3   +\cdots
,\\
\bar{a}_{33}(V)=& \frac{1}{8}I +\frac{3}{11200}V +
\frac{29}{8064000}V^2 -\frac{3667}{17740800000}V^3 +\cdots .\\
\end{aligned}
\end{equation}
We denote this method by SERKN3s4(2), which  is proved to satisfy
all the order conditions and symplectic conditions.

\section{Stability and Numerical experiments}\label{four}
\subsection{Stability regions}

In this section, we are concerned with the stability of the
 ERKN methods. This has been analyzed in
\cite{wu2012-2cpc} and  thence we just briefly recall here the
definitions. Consider the revised test equation:
\begin{equation}
y^{\prime\prime}(t)+\omega^{2}y(t)=-\epsilon y(t)\qquad
\mathrm{with}
\qquad \omega^{2}+ \epsilon>0, \label{ekrnP}%
\end{equation}
where $\omega$ represents an estimation of the dominant frequency
$\lambda$ and $\epsilon=\lambda^{2}-\omega^{2}$ is the error of that
estimation. Applying a multi-frequency ERKN method to \eqref{ekrnP}
produces
\[
\left(
\begin{array}
[c]{c}%
q_{n+1}\\
hq_{n+1}^{\prime}%
\end{array}
\right)  =S(V,z)\left(
\begin{array}
[c]{c}%
q_{n}\\
hq_{n}^{\prime}%
\end{array}
\right)  ,
\]
where the stability matrix $S(V,z)$ is given by
\[
S(V,z)=\left(
\begin{array}
[c]{cc}%
\phi_{0}(V)-z\bar{b}^{T}(V)N^{-1}\phi_{0}(c^{2}V) & \phi_{1}(V)\!-\!z\bar{b}%
^{T}(V)N^{-1}\big(c\cdot\phi_{1}(c^{2}V)\big)\\
-V\phi_{1}(V)\!-\!zb^{T}(V)N^{-1}\phi_{0}(c^{2}V) & \phi_{0}(V)\!-\!zb^{T}%
(V)N^{-1}\big(c\cdot\phi_{1}(c^{2}V)\big)
\end{array}
\right)
\]
with $V=h^{2}\omega^{2},\ z=h^{2}\epsilon $ and $N=I+z\bar{A}(V)$.
\begin{defi}
\label{sta properties} (See \cite{wu2012-2cpc}.) $R_{s}=\{(V,z)|\
V>0\ \textmd{and}\ \rho(S)<1\}$ is called the stability region of a
multi-frequency
 ERKN method and $R_{p}=\{(V,z)|\ V>0,\ \rho(S)=1\ \textmd{and}\
\mathrm{tr}(S)^{2}<4\det(S)\}$ is called the periodicity region of a
multi-frequency    ERKN method.
\end{defi}

\begin{figure}[ptbh]
\centering\tabcolsep=1mm
\begin{tabular}
[c]{ccc}%
\includegraphics[width=5cm]
{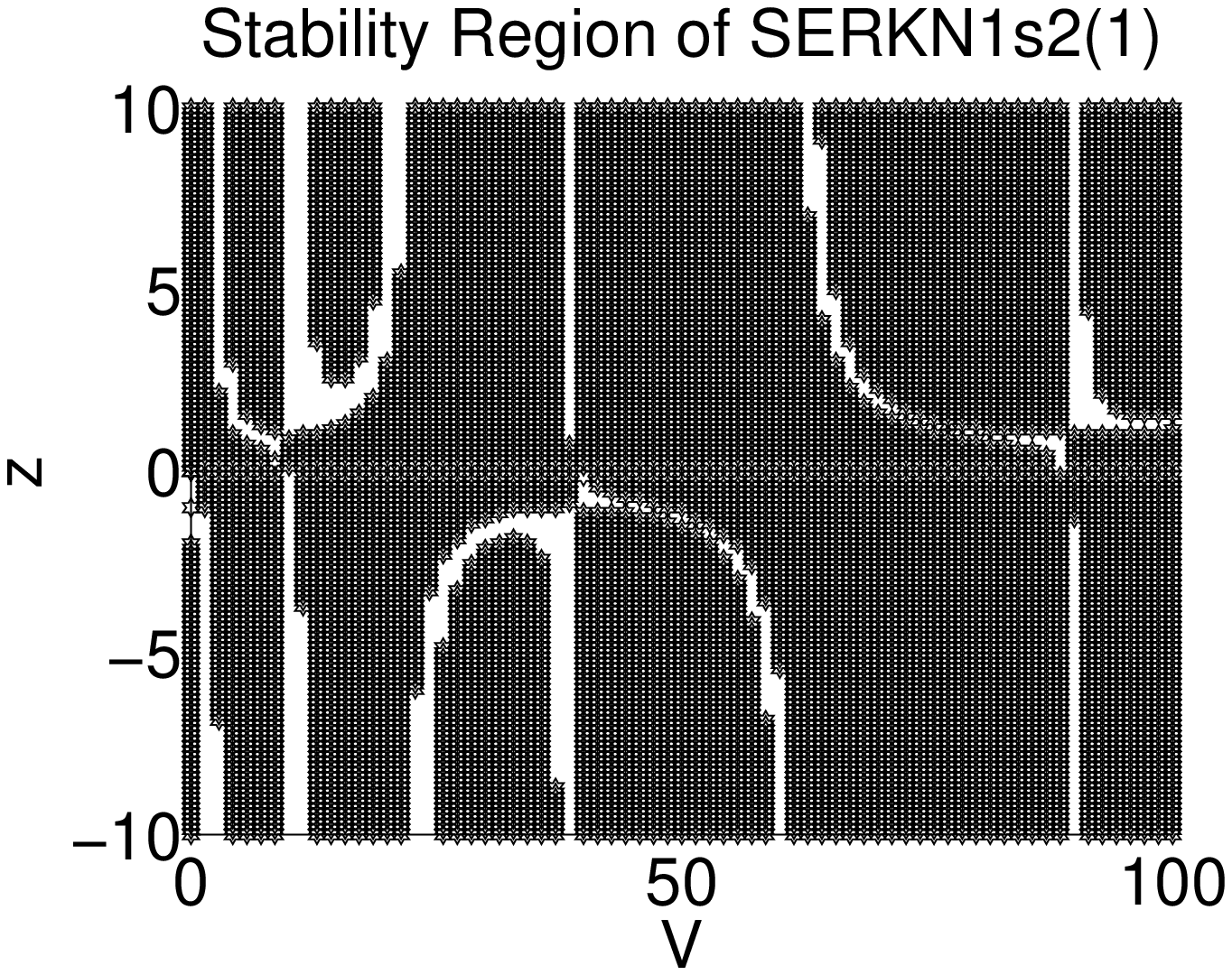} & \includegraphics[width=5cm]
{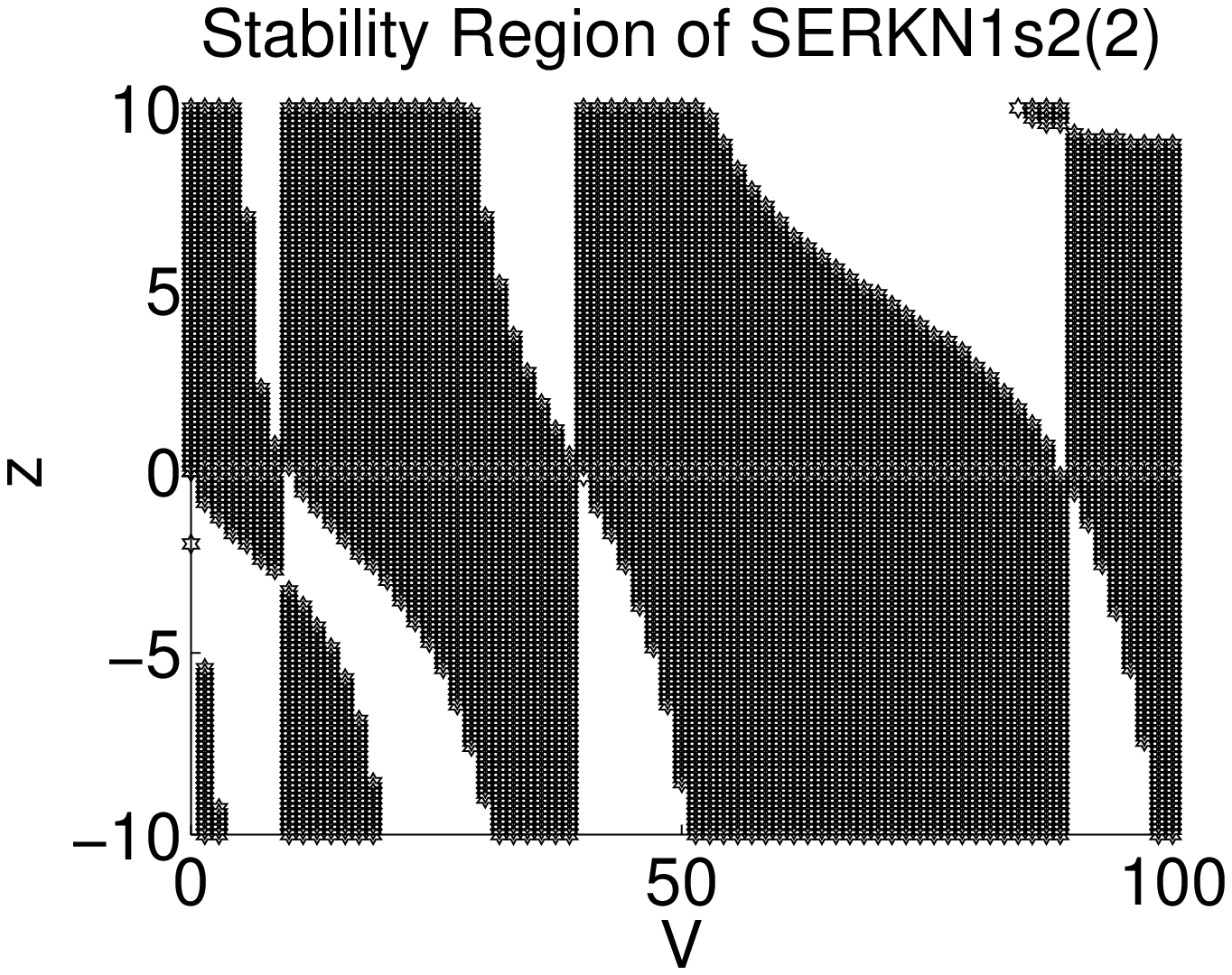} & \includegraphics[width=5cm]
{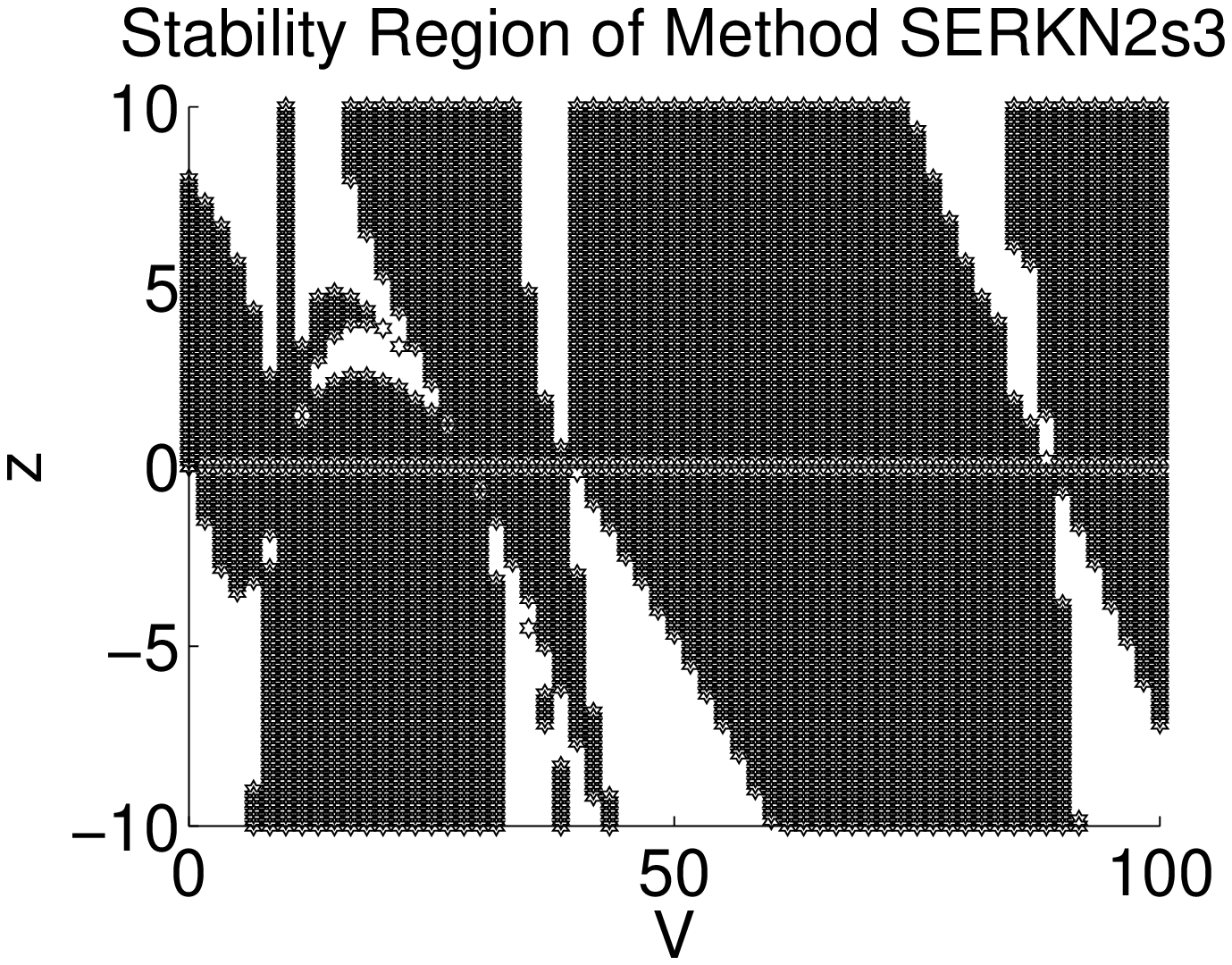}\\
\includegraphics[width=5cm]
{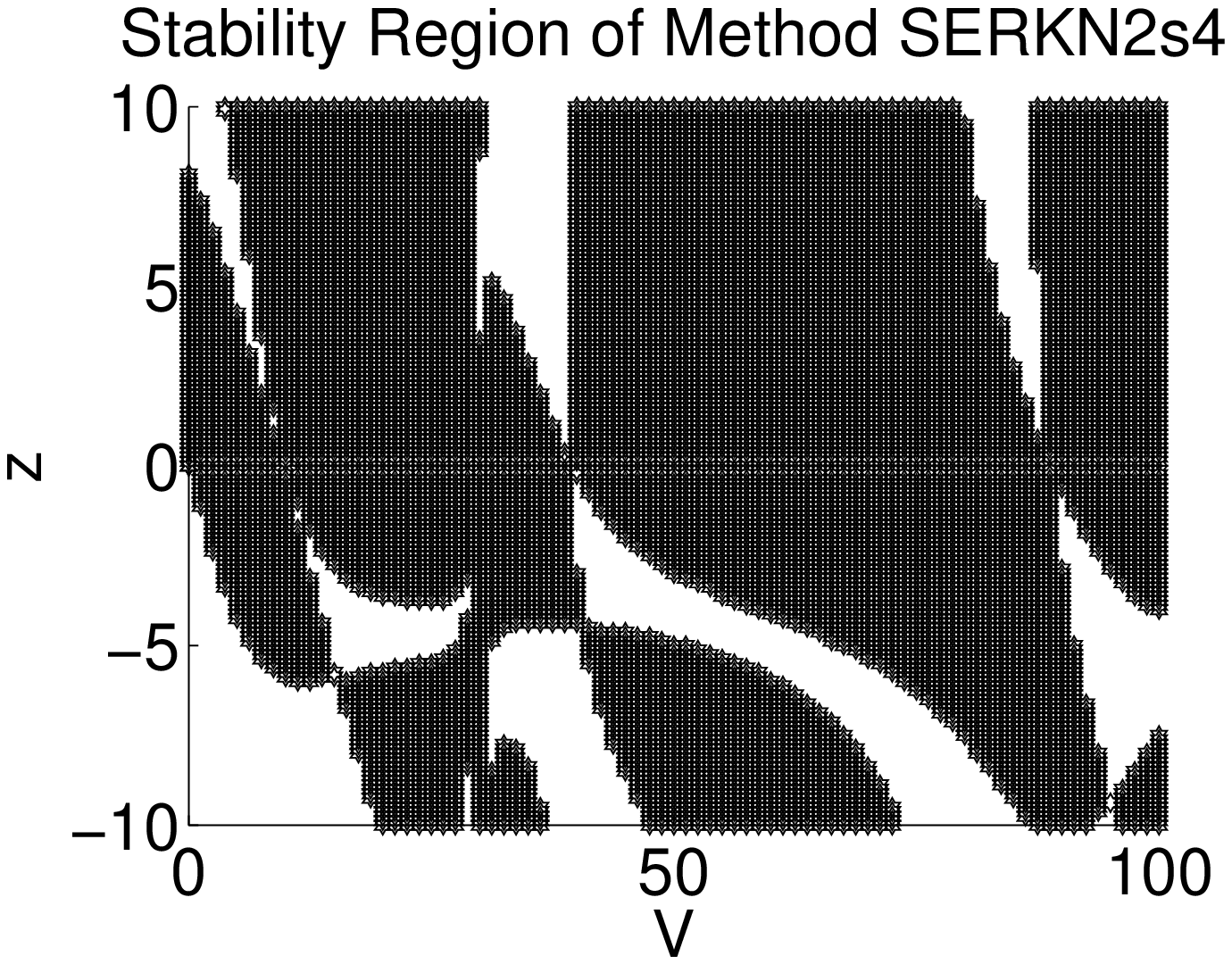} & \includegraphics[width=5cm]
{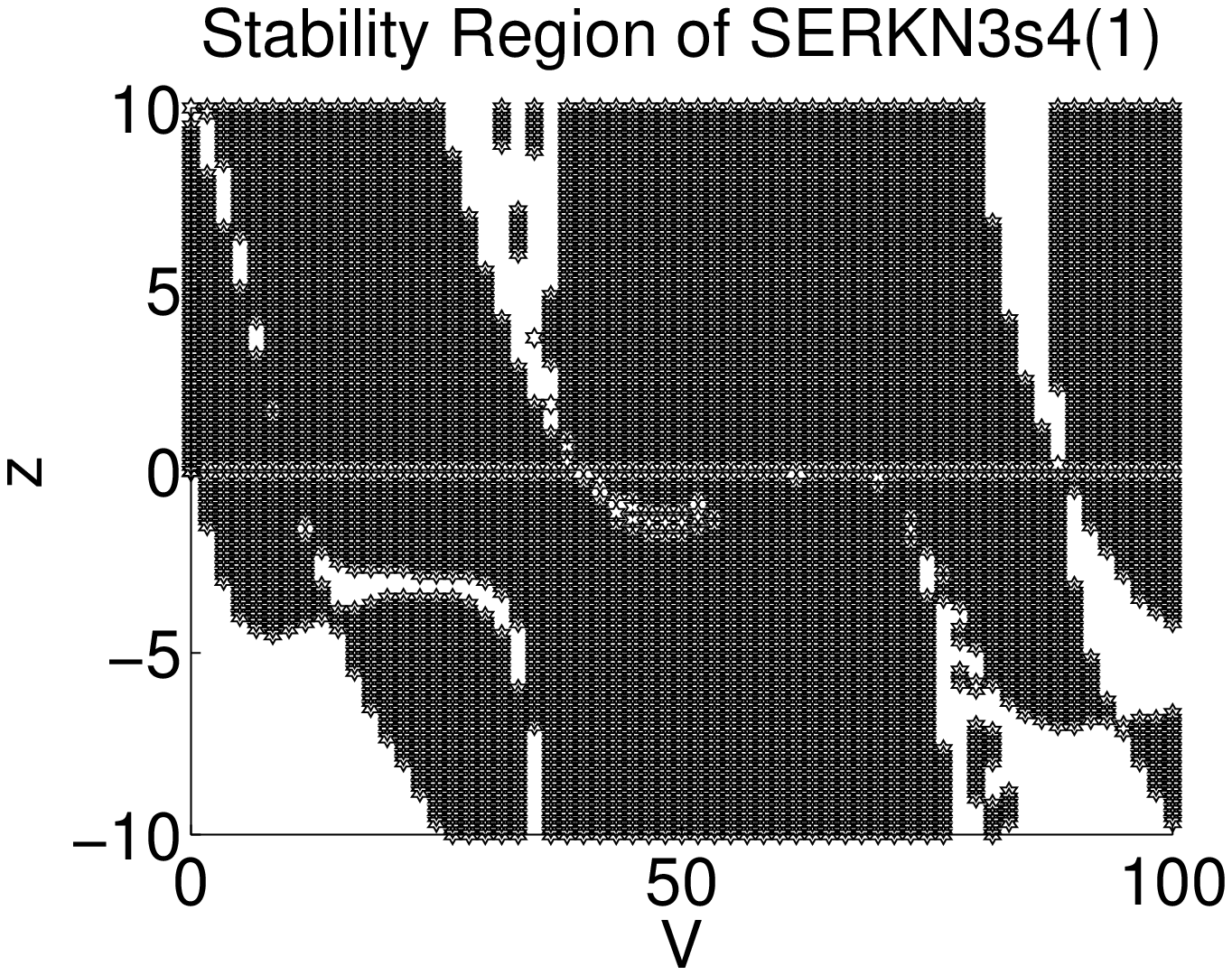} & \includegraphics[width=5cm]
{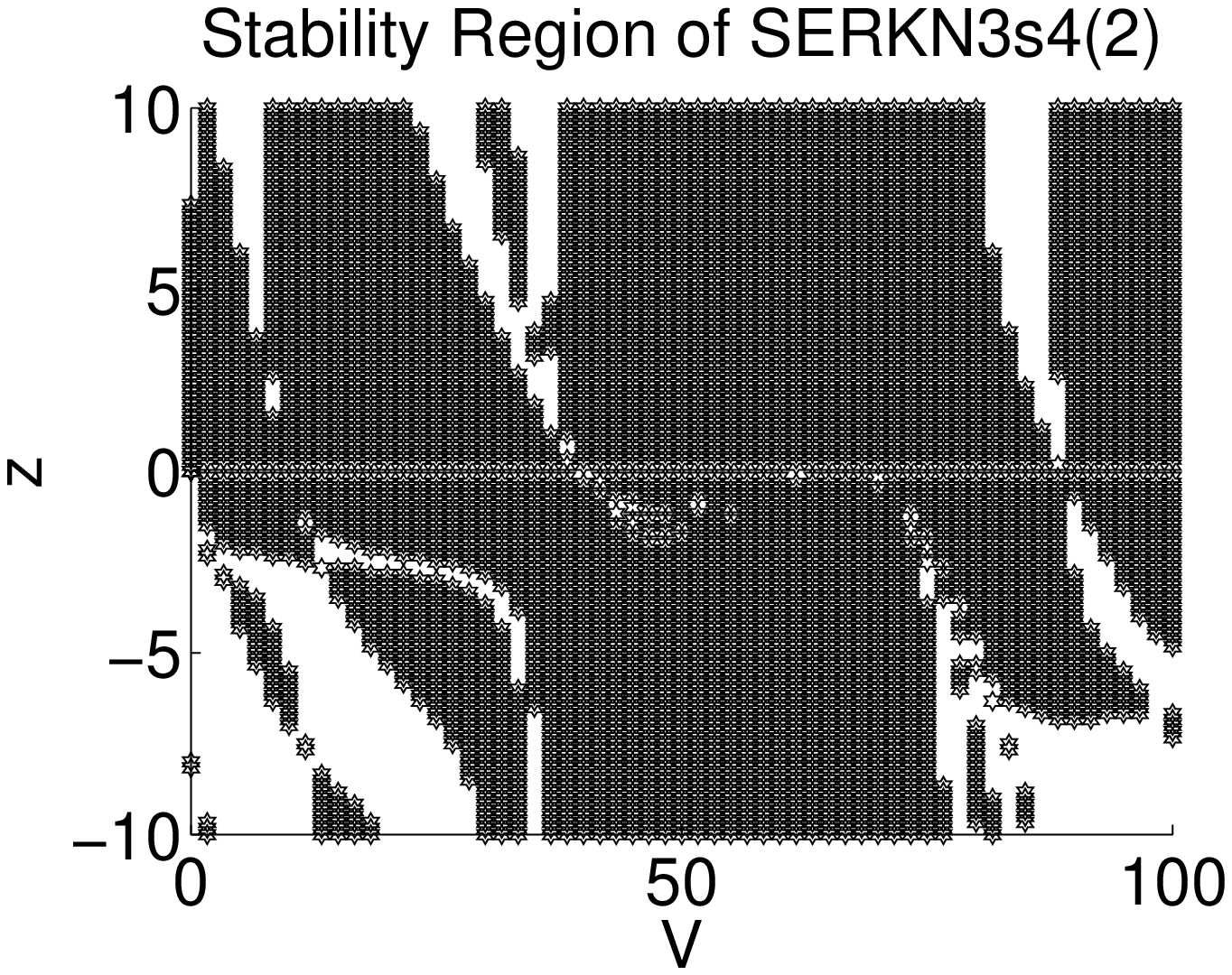}
\end{tabular}
\caption{{\protect\small Stability regions (shaded regions) for the obtained methods.}}%
\label{fig:mesrkn}%
\end{figure}

The stability regions of our methods are depicted in Figure
$\ref{fig:mesrkn}$.

\subsection{Numerical experiments}
In this section, in order to show the efficiency of our new methods
compared with their corresponding methods, we give three numerical
experiments and their results. The methods for comparisons are:
\begin{itemize}\itemsep=-0.2mm
\item SERKN1s2: the one-stage  diagonal implicit symplectic ERKN method (SERKN1s2(1)) of order two derived in Subsection\ref{ghhgg};

\item SERKN2s3: the two-stage  diagonal implicit symplectic ERKN method of order three derived in Subsection \ref{112233};

\item SERKN3s4: the three-stage diagonal implicit symplectic ERKN method (SERKN3s4(1)) of order four derived in Subsection \ref{334455};

\item RKN1s2: the one-stage  diagonal implicit symplectic RKN method of order two obtained by
 letting $V \rightarrow \textbf{0}$  of the method SERKN1s2;
\item RKN2s3: the two-stage  diagonal implicit symplectic RKN method of order three obtained by
 letting $V \rightarrow \textbf{0}$  of the method SERKN2s3;
\item RKN3s4: the three-stage  diagonal implicit symplectic RKN method of order four obtained by
 letting $V \rightarrow \textbf{0}$  of the method SERKN3s4.
\end{itemize}
The numerical experiments have been carried out on a personal
computer and the algorithm has been implemented by using the
MATLAB-R2010a.

\begin{figure}[ptbh]
\centering\tabcolsep=1mm
\begin{tabular}
[c]{ccc}%
\includegraphics[width=5cm,height=5.5cm]{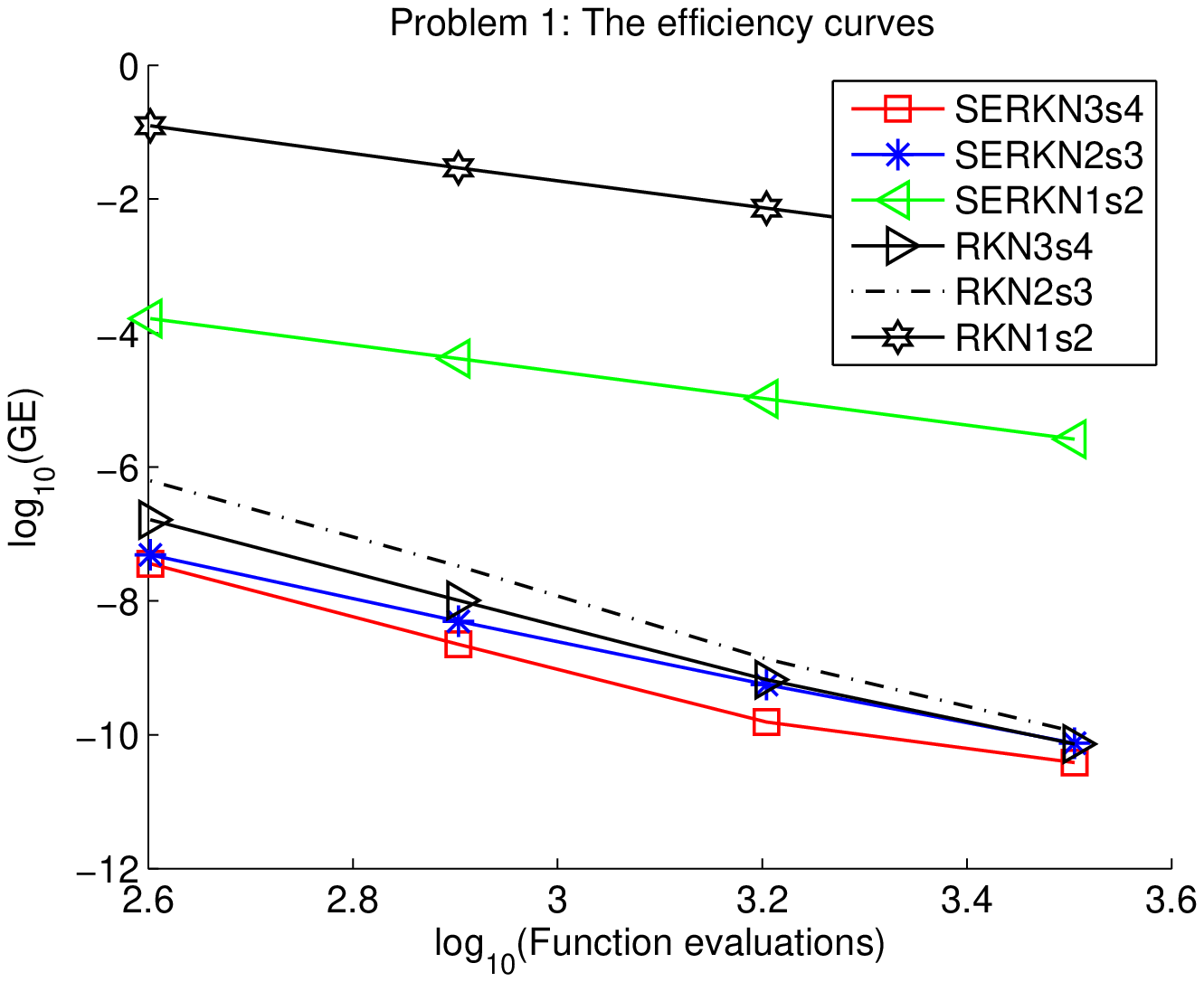} &\includegraphics[width=5cm,height=5.5cm]{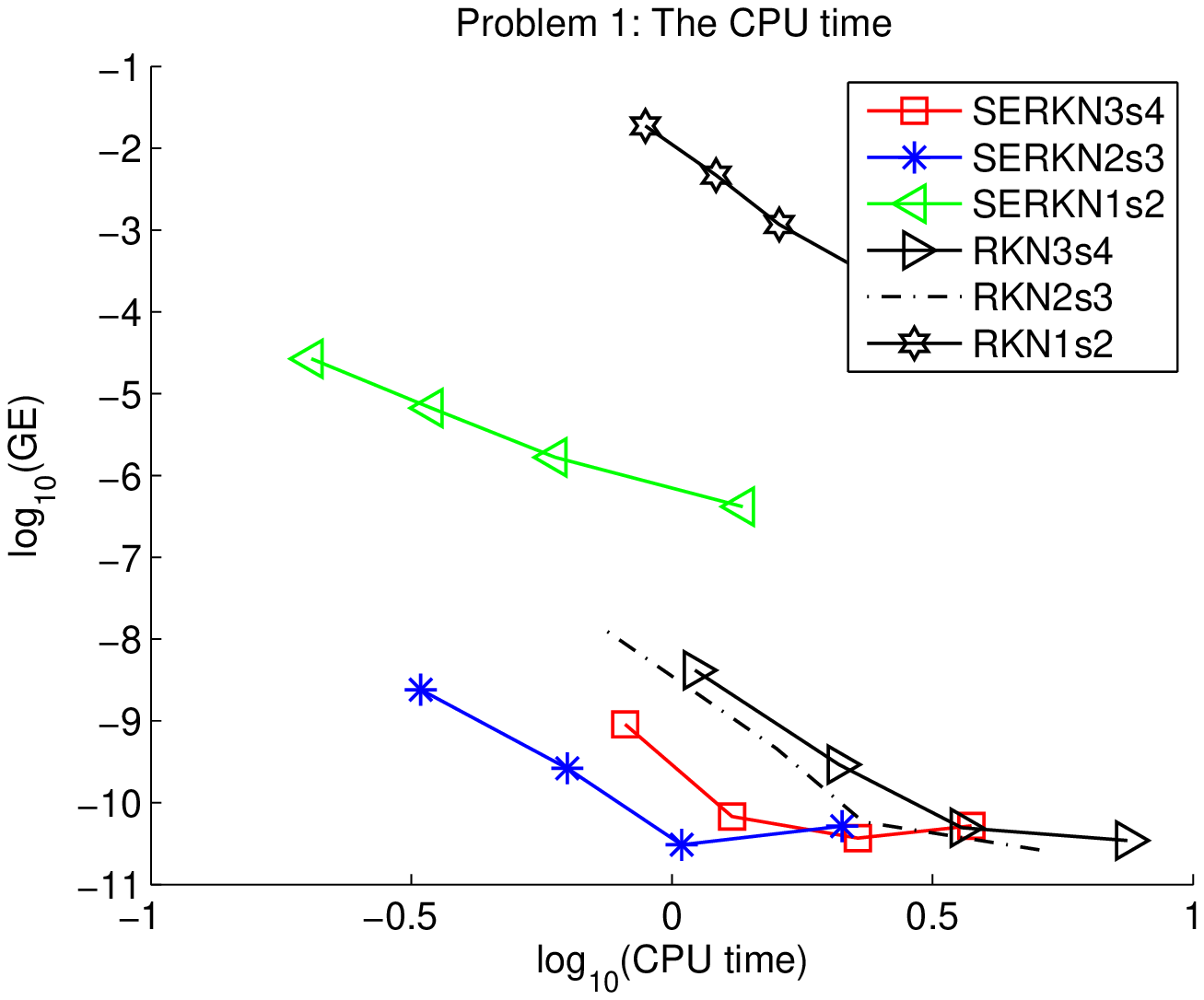} &\includegraphics[width=5cm,height=5.5cm]{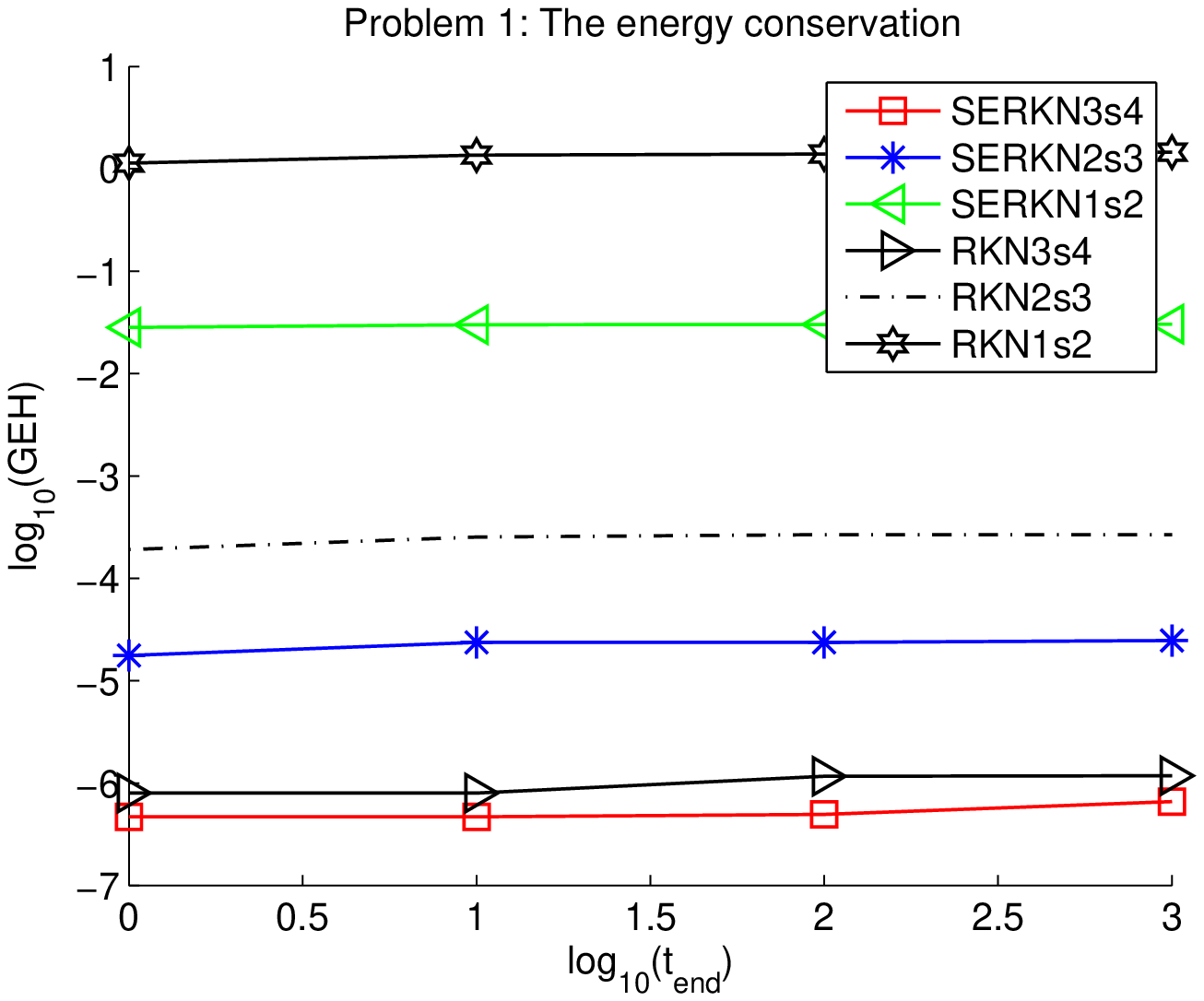}\\
{\small (i)} & {\small (ii)} & {\small (iii)}%
\end{tabular}
\caption{Results for Problem 1. (i): The logarithm of the global
error ($GE$) over the integration interval against the logarithm of
the number of function evaluations.  (ii): The logarithm of the
global error ($GE$) over the integration interval against the CPU
time. (iii):\ The logarithm of the maximum global error of
Hamiltonian
$GEH=\max|H_{n}-H_{0}|$ against $\log_{10}(t_{\mathrm{end}})$.}%
\label{fig:problem1}%
\end{figure}

\vskip1mm \noindent   \textbf{Problem 1.} Consider the sine-Gordon
 equation with periodic boundary conditions (see \cite{smmmeijer1987})
\begin{eqnarray*}
\begin{array}{ll}
\dfrac{\partial^2u}{\partial t^2}=\dfrac{\partial^2u}{\partial
x^2}-\sin u,\ \ \ \ \ -1<x<1,\ \ \ \ \ t>0,\ \ \
u(-1,t)=u(1,t).\\[0.3cm]

\end{array}
\end{eqnarray*}
 We carry out a semi-discretization on the spatial by using
 second-order symmetric differences and obtain the following system
 of second-order ODEs in time
\begin{eqnarray*}
\begin{array}{ll}
\dfrac{d^2U}{dt^2}+MU=F(t,U),\ \ \ 0<t\leq t_{\mathrm{end}},
\end{array}\label{pro4}
\end{eqnarray*}
where $U(t)=\big(u_{1}(t),\ldots,u_{N}(t)\big)^{T}$ with
$u_{i}(t)\approx u(x_{i},t)$, $i=1,\ldots,N,$
\begin{eqnarray*}
M=\dfrac{1}{\Delta x^2}\left(
\begin{array}
[c]{ccccc}%
2&-1 && &-1\\
-1 &2 & -1&  &  \\
 &\ddots&\ddots&\ddots& \\
&&-1 &2 & -1\\
-1& & &-1&2 \\
\end{array}
\right)
\end{eqnarray*}
with $\Delta x= 1/N$,and  $x_i =-1+ i\Delta x$,and
$F(t,U)=-\sin(U)=-\big(u_{1},\ldots,u_{N}\big)^{T}.$ The
  Hamiltonian of this system is
\begin{equation*}
H(U',U)=\dfrac{1}{2}U'^{T}U'+\dfrac{1}{2}U^{T}MU-\cos(u_{1})-\cos(u_{2})-\ldots-\cos(u_{N}).
\end{equation*}
We take the initial conditions as
\[
U(0)=(\pi)_{i=1}^{N},\ \ \ U_{t}(0)=\sqrt{N}\Big(0.01+\sin(\dfrac{2 \pi i}%
{N})\Big)_{i=1}^{N}
\]
with $N=32$. We integrate this problem in the interval $[0,10]$ with
stepsizes $h=1/(20\times 2^{i}), \ \ i=1,2,3,4$. Figure
\ref{fig:problem1} (i) shows the global errors. We then solve this
problem in the interval $[0,10]$ with the stepsize $h=1/(100\times
2^{i}),\ \ i=1,2,3,4$ and show the global errors against the CPU
time in Figure \ref{fig:problem1} (ii). Finally we integrate this
problem with a fixed stepsize $h=1/40$ in the interval
$[0,t_{end}],\ \ t_{end}=10^{i}$ with $i=0,1,2,3$. The results of
energy conservation are presented in Figure \ref{fig:problem1}
(iii).

\begin{figure}[ptbh]
\centering\tabcolsep=1mm
\begin{tabular}
[c]{ccc}%
\includegraphics[width=5cm,height=5.5cm]{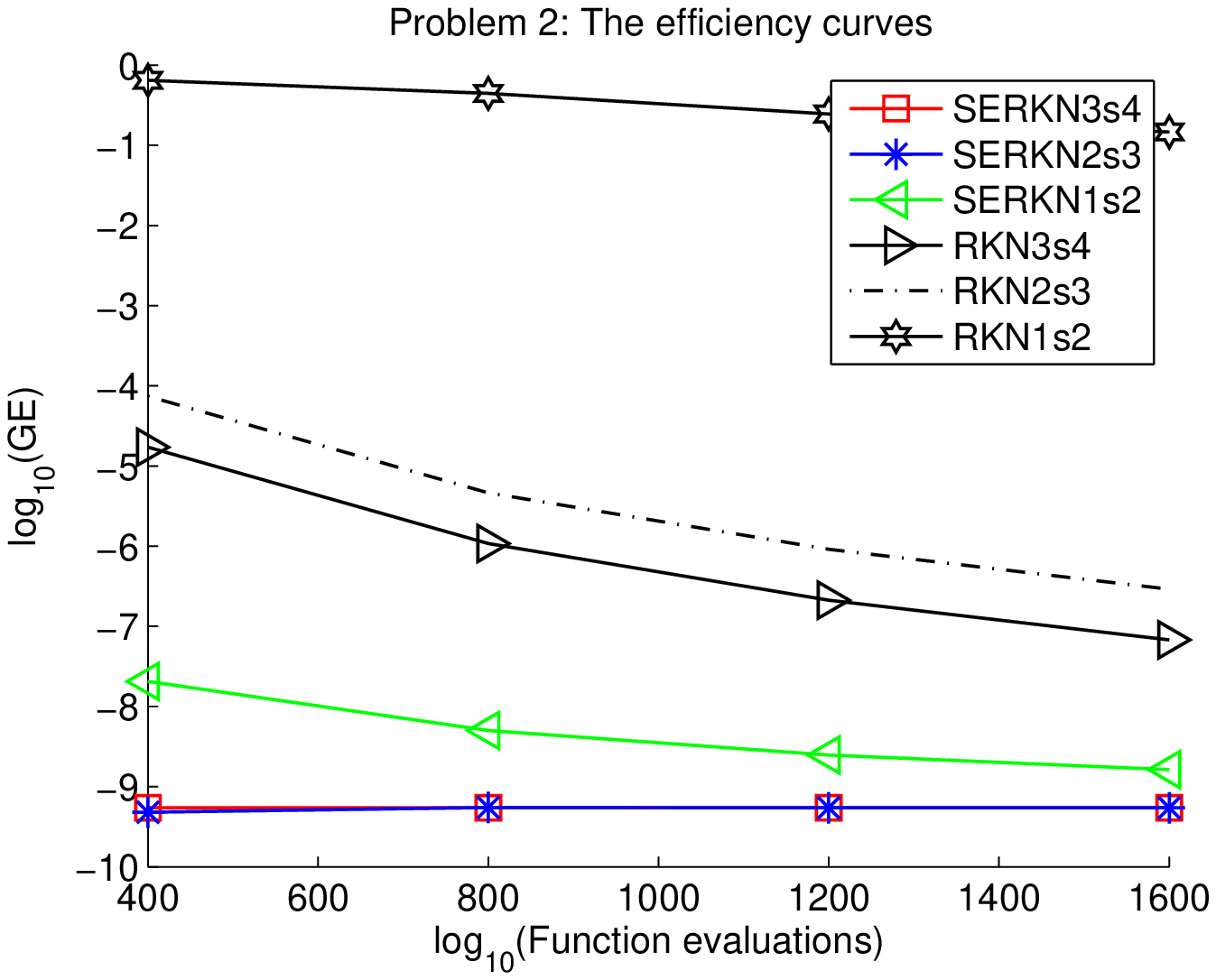} &\includegraphics[width=5cm,height=5.5cm]{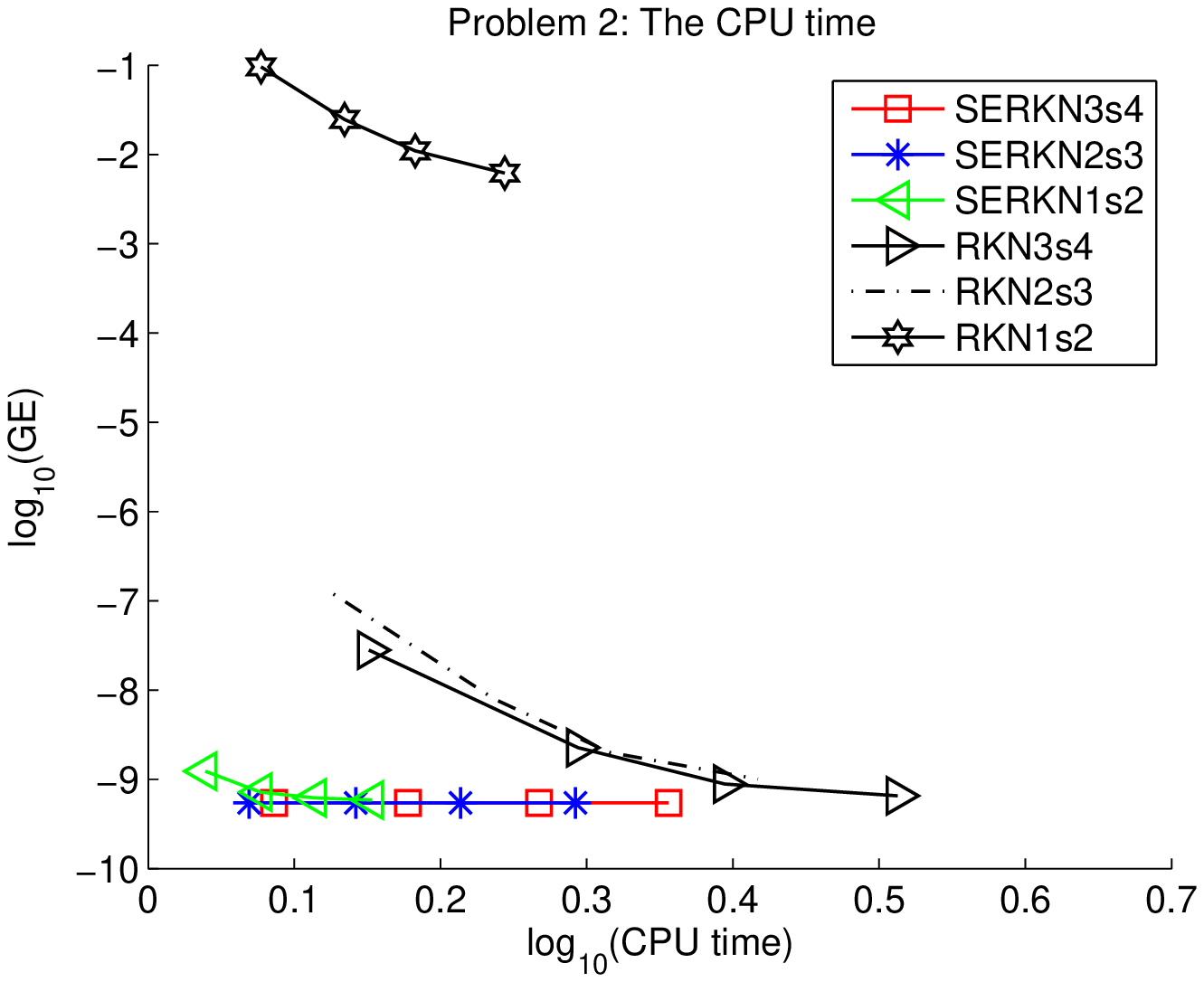} &\includegraphics[width=5cm,height=5.5cm]{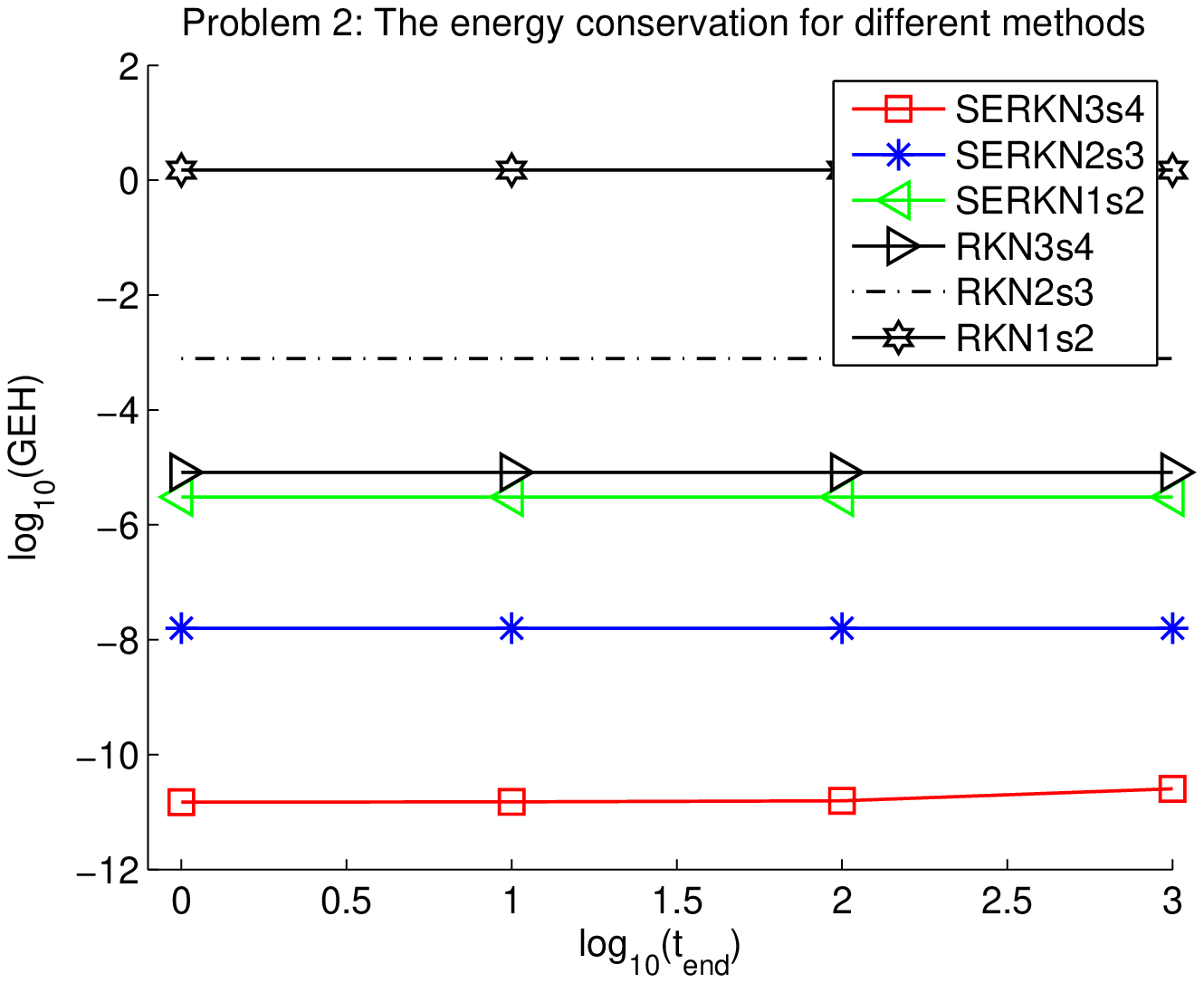}\\
{\small (i)} & {\small (ii)} & {\small (iii)}%
\end{tabular}
\caption{Results for Problem 2. (i): The logarithm of the global
error ($GE$) over the integration interval against the logarithm of
the number of function evaluations. (ii): The logarithm of the
global error ($GE$) over the integration interval against the CPU
time. (iii):\ The logarithm of the maximum global error of
Hamiltonian
$GEH=\max|H_{n}-H_{0}|$ against $\log_{10}(t_{\mathrm{end}})$.}%
\label{fig:problem2}%
\end{figure}

\noindent\vskip3mm \noindent\textbf{Problem 2.} Consider the Duffing
equation
\begin{equation*}
q''+100q=k^{2}\big(2q^{3}-q\big),  \qquad
q(0)=0,\ \ q'(0)=10,\qquad t\in[0,t_{\mathrm{end}}],\label{prob}%
\end{equation*}
with $0\leqslant k<10.$  The
  Hamiltonian of this system is
\begin{equation*}
H(q,q')=\frac{1}{2}q'^{2}+50q^{2}-k^{2}(\frac{1}{2}q^{4}-\frac{1}{2}q^{2}).
\end{equation*}
The analytic solution of this initial value problem is given by
$q(t)=sn(10t,k/10)$, and represent a periodic motion in terms of the
Jacobian elliptic function $sn$. In this test we choose the
parameter values $k=0.03,$ and integrate this problem with the
stepsize $h=1/(200\times i),\ \ i=1,2,3,4$ in the interval $[0,10]$.
See Figure \ref{fig:problem2} (i) for the  efficiency curves. We
then solve this problem in the interval $[0,10]$ with the stepsize
$h=1/(40\times i),\ \ i=1,2,3,4$ and show the global errors against
the CPU time in Figure \ref{fig:problem2} (ii). Finally we integrate
this problem with a fixed stepsize $h=1/50$ in the interval
$[0,t_{\mathrm{end}}],\ \ t_{\mathrm{end}}=10^{i}$ with $i=0,1,2,3$.
The results of energy conservation are presented in Figure.
\ref{fig:problem2} (iii).

\begin{figure}[ptbh]
\centering\tabcolsep=1mm
\begin{tabular}
[c]{ccc}%
\includegraphics[width=5cm,height=5.5cm]{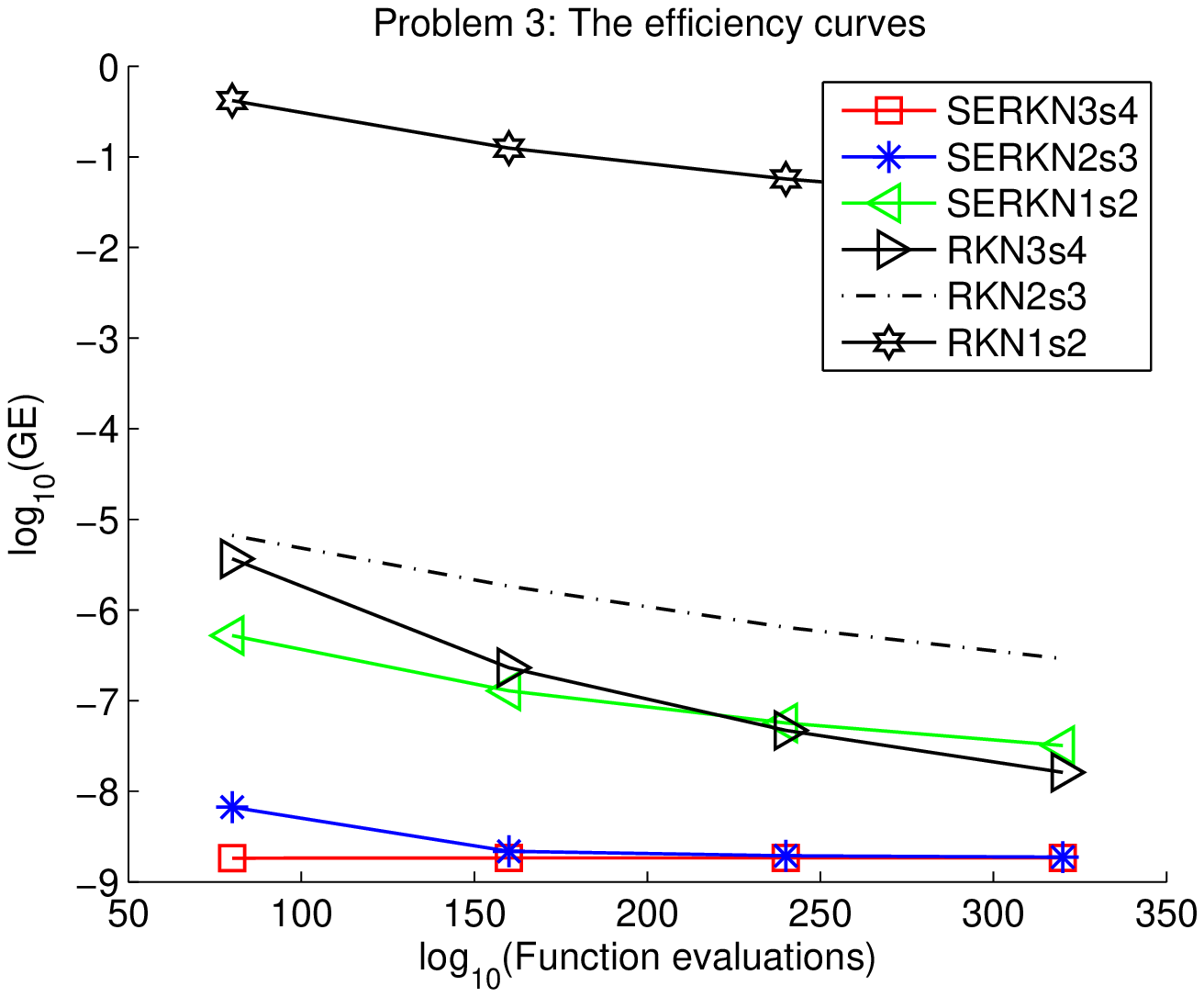} & \includegraphics[width=5cm,height=5.5cm]{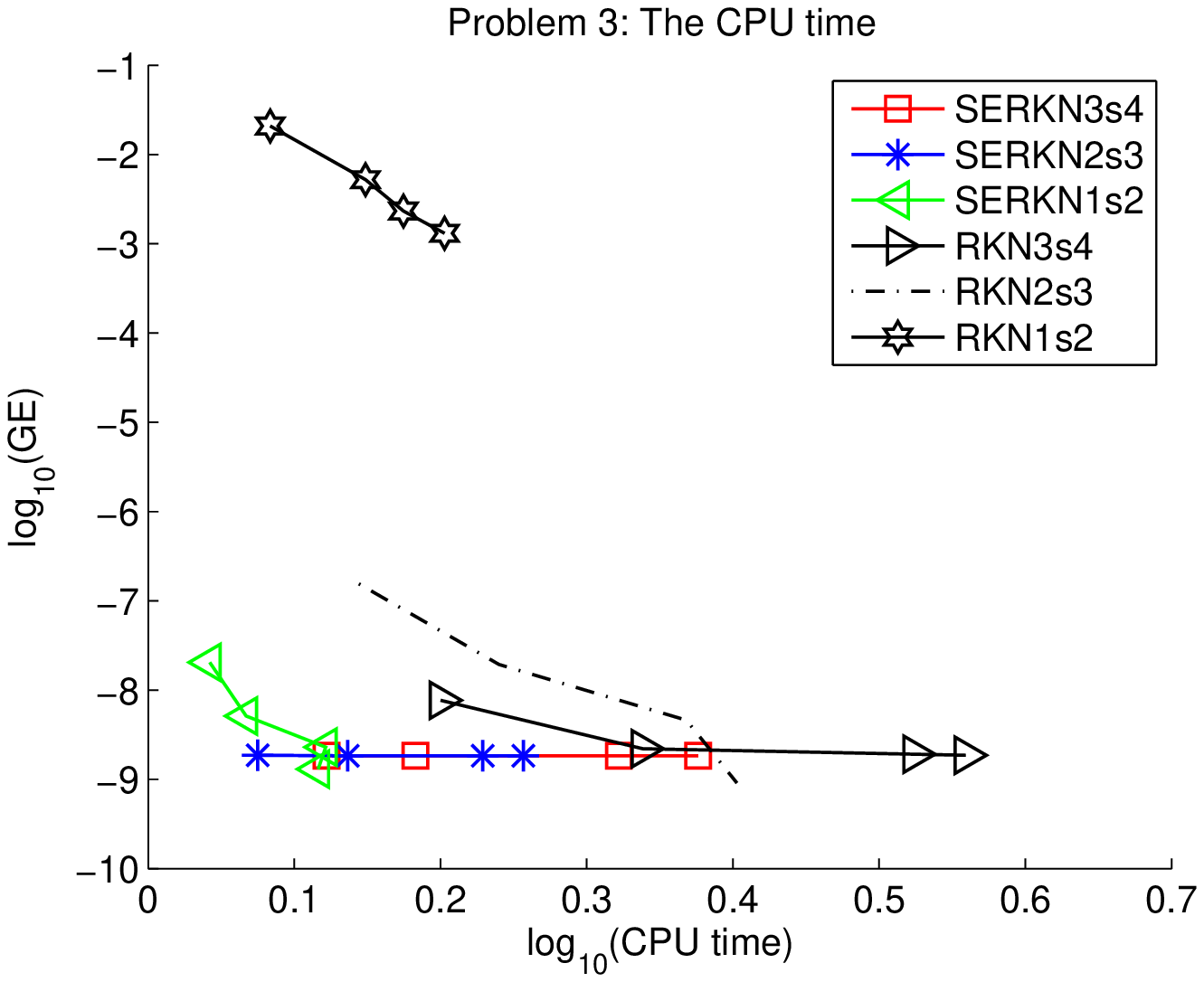} & \includegraphics[width=5cm,height=5.5cm]{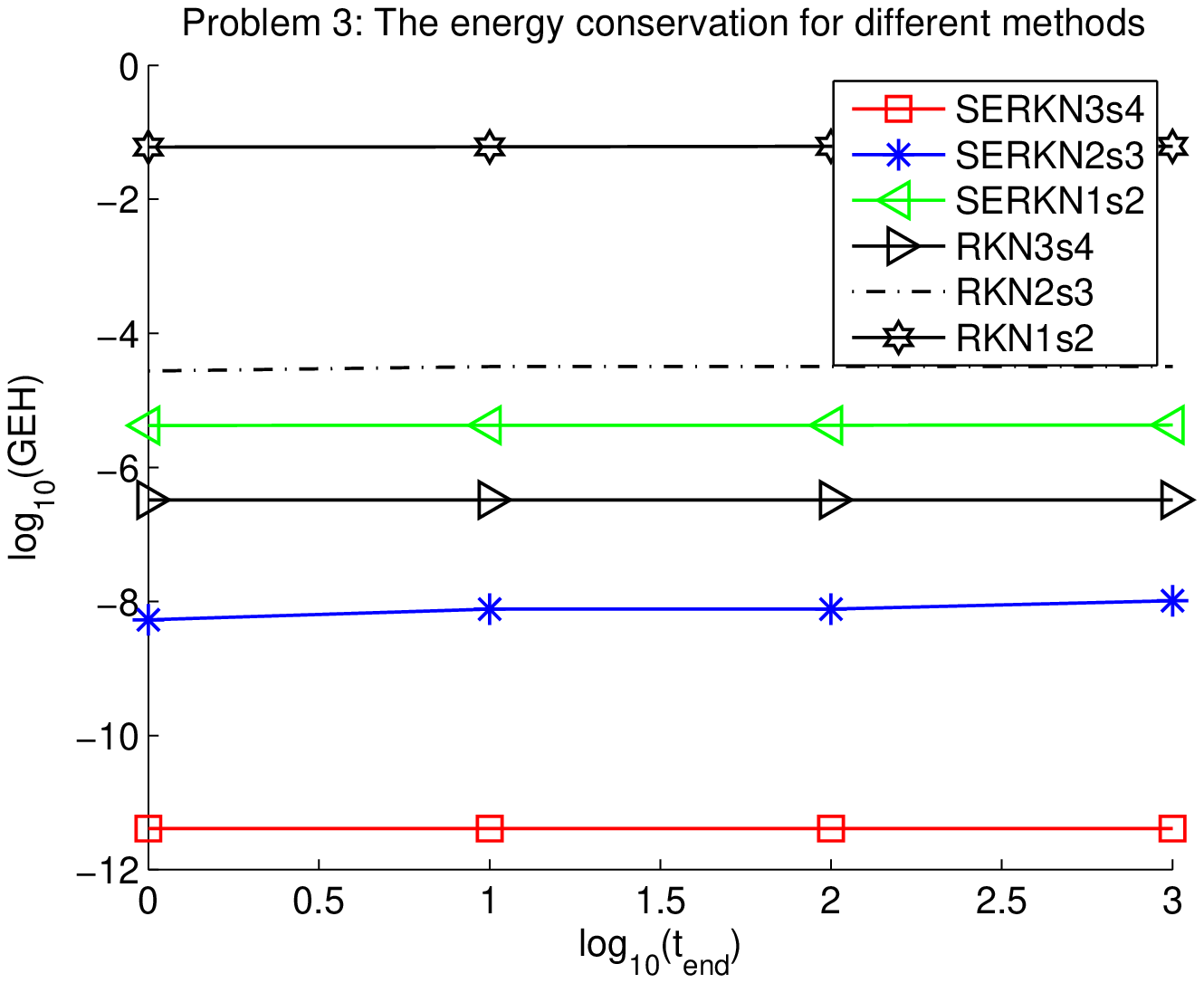}\\
{\small (i)} & {\small (ii)} & {\small (iii)}%
\end{tabular}
\caption{Results for Problem 3. (i): The logarithm of the global
error ($GE$) over the integration interval against the logarithm of
the number of function evaluations. (ii): The logarithm of the
global error ($GE$) over the integration interval against the CPU
time. (iii):\ The logarithm of the maximum global error of
Hamiltonian
$GEH=\max|H_{n}-H_{0}|$ against $\log_{10}(t_{\mathrm{end}})$.}%
\label{fig:problem3}%
\end{figure}

\noindent\vskip3mm \noindent\textbf{Problem 3.} Consider the model
for stellar orbits in a galaxy (see
\cite{Kevorkian1981,Kevorkian1996})
\begin{equation*}\begin{array}
[c]{ll} q''_{1}(t)+a^{2}q_{1}(t)=\epsilon q_{2}^{2}(t),\qquad\qquad\
q_{1}(0)=1,\ \
q'_{1}(0)=0,\\
q''_{2}(t)+b^{2}q_{2}(t)=2\epsilon q_{1}(t)q_{2}(t),\qquad
q_{2}(0)=1,\ \
q'_{2}(0)=0,\\
\end{array}
\label{prooo}%
\end{equation*}
where $q_{1}$ stands for the radial displacement of the orbit of a
star from a reference circular orbit, and $q_{2}$ stands for the
deviation of the orbit from the galactic plane. The time variable
$t$ actually denotes the angle of the planets in a reference
coordinate system. We choose $a=2,$ $b=1.$ The Hamiltonian of this
system is
\begin{equation*}
H(q,q')=\dfrac{1}{2}(q_{1}'^{2}+q_{2}'^{2})+\dfrac{1}{2}(4q_{1}^2+q_{2}^2)-\epsilon
q_{1}q_{2}^{2}.
\end{equation*}
The problem has been solved on the interval $[0,1000]$ with
$\epsilon=10^{-3}.$    We integrate this problem in the interval
$[0,10]$ with stepsizes $h=1/(8\times i), \ \ i=1,2,3,4$ and
efficiency curves  are presented in Figure \ref{fig:problem3}(i). At
the meantime, we solve this problem in the interval $[0,10]$ with
the stepsize $h=1/(40\times i),\ \ i=1,2,3,4$ and show the global
errors against the CPU time in Figure \ref{fig:problem3} (ii).
Finally we integrate this problem with a fixed stepsize $h=1/10$ in
the interval $[0,t_{\mathrm{end}}],\ \ t_{\mathrm{end}}=10^{i}$ with
$i=0,1,2,3$. The results of energy conservation are presented in
Figure \ref{fig:problem3} (iii).

It follows from the numerical results that our novel methods are
very promising as compared with their corresponding RKN methods.
\section{Conclusions}\label{five}
In this paper,  based on symplecticity conditions and order
conditions, we obtain one-stage of order two, two-stage of order
three and three-stage of order  four diagonal implicit symplectic
ERKN methods
  for solving the oscillatory Hamiltonian system
$\eqref{a}$. We also discuss the stability of the new methods.
Furthermore, numerical experiments are performed in comparsion with
their corresponding RKN methods in the scientific literature. The
remarkable efficiency of the new methods are shown by the numerical
results.

\end{document}